\documentclass[10pt]{article}

\usepackage[backend=biber, style=alphabetic, citestyle=alphabetic, backref=true]{biblatex}
\DeclareSourcemap{
  \maps[datatype=bibtex]{
    \map[overwrite]{
      \step[fieldsource=doi, final]
      \step[fieldset=url, null]
      \step[fieldset=eprint, null]
    }  
  }
}

\usepackage{booktabs}
\usepackage[margin=1in]{geometry}

\usepackage[usenames,dvipsnames,svgnames,table]{xcolor}
\usepackage{graphicx, float, wrapfig, caption, subcaption}
\usepackage{pdflscape}
\graphicspath{
{"/Users/rynebeeson/Box/LaTeX Images/"}
{"/Users/rynebeeson/Box Sync/LaTeX Images/"}
}

\usepackage[normalem]{ulem}

\usepackage{enumerate}

\usepackage{yfonts}

\usepackage[mathscr]{eucal}

\usepackage[final]{hyperref}
\hypersetup{colorlinks=true}
\hypersetup{linkcolor=Emerald}
\hypersetup{citecolor=cyan}
\hypersetup{menucolor=red}

\usepackage{beesonmacros}
\newcommand{\eps}{\epsilon}
\newcommand{\XOe}{X^{0, \epsilon}}
\newcommand{\Xe}{X^\epsilon}
\newcommand{\Ze}{Z^\epsilon}
\newcommand{\Zex}{Z^{\epsilon, x}}
\newcommand{\Ye}{Y^\epsilon}

\newcommand{\mYe}{\mathcal{Y}^\epsilon}
\newcommand{\mF}{\mathcal{F}}
\newcommand{\wt}[1]{\widetilde{#1}}
\newcommand{\Pe}{\mathbb{P}^\epsilon}
\newcommand{\mG}{\mathcal{G}}

\newcommand{\wDe}{\wt{D}^\epsilon}

\newcommand{\Qo}{\mathbb{Q}^0}

\newcommand{\B}[1][t]{\overleftarrow{B}_{#1}}

\newcommand{\Zextz}[1][s]{Z^{\epsilon, x; (t, z)}_{#1}}
\newcommand{\Xetx}[1][s]{X^{\epsilon; (t, x)}_{#1}}
\newcommand{\Zetz}[1][s]{Z^{\epsilon; (t, z)}_{#1}}
\newcommand{\F}[1][t, T]{\mF_{#1}^B}

\renewcommand{\norm}[2]{| #1 |_{#2}}

\usepackage{fancyhdr}
\usepackage[yyyymmdd,hhmmss]{datetime}
\newcommand{\orcid}[1]{\href{https://orcid.org/#1}{\includegraphics[scale=.014]{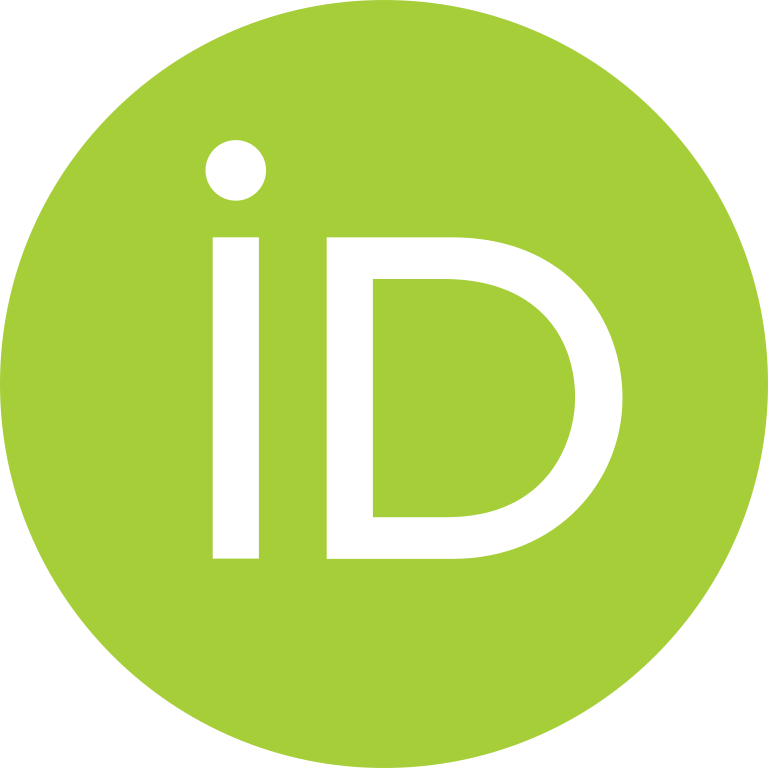}}}

\title{Quantitative Convergence of the Filter Solution for Multiple Timescale Nonlinear Systems with Coarse-Grain Correlated Noise}
\date{\today}
\author{Ryne Beeson\footnote{Princeton University} \orcid{0000-0003-2176-0976},
N. Sri Namachchivaya\footnote{University of Waterloo}, 
and Nicolas Perkowski\footnote{Freie Universit\"at Berlin}}

\begin{document}

\maketitle

\begin{abstract}
In this paper we prove a rate of convergence for the continuous time filtering solution of a multiple timescale correlated nonlinear system to a lower dimensional filtering equation in the limit of large timescale separation. 
Correlation is assumed to occur between the slow signal and observation processes. 
Convergence is almost sure in the weak topology.  
An asymptotic expansion of the dual process for the solution to the Zakai equation, and probabilistic representation using backward doubly stochastic differential equations is leveraged to prove the result. 
\end{abstract}

\section{Introduction}

In this paper we prove a rate of convergence for the continuous time filtering solution of a multiple timescale and correlated nonlinear system to a lower dimensional filtering equation. 
The coupled system of stochastic differential equations (SDEs) that we consider is as follows, 
\EquationAligned{
\label{equation: multiscale signal process}
d\Xe_t 
&= b(\Xe_t, \Ze_t) dt + \sigma(\Xe_t, \Ze_t) dW_t, \\
d\Ze_t 
&= \frac{1}{\epsilon^2} f(\Xe_t, \Ze_t) dt + \frac{1}{\epsilon} g(\Xe_t, \Ze_t) dV_t.
}
We denote the infinitesimal generator of $(\Xe, \Ze)$ as $\mG^\epsilon$.
The process $(\Xe, \Ze)$ is known as the signal process and $\epsilon \in (0, 1)$ is a timescale parameter such that $\Ze$ is a fast process and $\Xe$ is a slow process. 
In filtering theory, we consider the signal process to be non-observable, and instead have indirect measurements of $(\Xe, \Ze)$ via the noisy observation process,
\begin{align*}
d\Ye_t = h(\Xe_t, \Ze_t) dt + \alpha dW_t + \gamma dU_t.
\end{align*}
We assume $W, V, U$ are independent Brownian motions, and the presence of $\R{d \times w} \ni \alpha \neq 0$ indicates correlation between the observation and slow (coarse-grain) process.
The goal in filtering theory is then to calculate the conditional distribution of $(\Xe, \Ze)$ given the observation history generated from $\Ye$, which we denote by $\pi^\epsilon$. 
At each time $t > 0$, $\pi^\epsilon_t$ is a random probability measure on the space $\R{m} \times \R{n}$ and acts on test functions $\varphi : \R{m} \times \R{n} \rightarrow \R{}$ by integration $\pi^\epsilon_t(\varphi) = \int \varphi(x, z) \pi^\epsilon_t(dx, dz)$. 

The question of this paper is then motivated by the well known result in homogenization of stochastic differential equations that if for every fixed $x$, the solution $Z^x$ of
\begin{align*}
dZ^x_t = f(x, Z^x_t) dt + g(x, Z^x_t) dV_t,
\end{align*}
is ergodic with stationary distribution $\mu_\infty(x)$, then under appropriate assumptions, the process $\Xe$ converges in distribution to a Markov process $X^0$ with infinitesimal generator $\overbar{\mG_S}$ in the limit as $\epsilon \rightarrow 0$ \cite{Papanicolaou:1977dtc, Pardoux:2003, Khasminskii:2005jde}.
Therefore, if we are only interested in statistics of $\Xe$ (i.e., estimation of test functions $\varphi : \R{m} \rightarrow \R{}$), then it would be computationally advantageous to know if $\pi^{\epsilon, x} \Rightarrow \pi^0$ converges weakly to a lower dimensional filtering equation; $\pi^0_t$ being a random probability measure for each time $t$ on $\R{m}$ and $\pi^{\epsilon, x}$ being the $x$-marginal of $\pi^\epsilon$. 

Filtering theory has widespread applications in many fields including various disciplines of engineering for decision and control systems, the geosciences, weather and climate prediction. 
In many of these fields, it is not uncommon to have physics based models with multiple timescales as seen in Eq. \ref{equation: multiscale signal process}, and also have the case were estimation of the slow process is solely of interest; for example the estimation of the ocean temperature, which is necessary for climate prediction, but the ocean model may also be coupled to a fast atmospheric model. 
Knowing that mathematically $\pi^{\epsilon, x} \Rightarrow \pi^0$ in the limit as $\epsilon \rightarrow 0$, enables practitioners to devise more efficient methods for estimation of the slow process without great loss of accuracy (see for instance \cite{Park:2011jam, Kang:2012hj, Berry:2014, Yeong:2020}).

There are several papers providing results for $\pi^{\epsilon, x} \rightarrow \pi^0$ (or the associated unnormalized conditional measure or density versions) on variations of the multiple timescale filtering problem. 
In \cite{Park:2010n}, $(\Xe, \Ze)$ is a two dimensional process with no drift in the fast component, no intermediate scale, and no correlation. 
The authors made use of a representation of the slow component by a time-changed Brownian motion under a suitable measure to yield weak convergence of the filter. 
Homogenization of the nonlinear filter was studied in \cite{Bensoussan:1986s} and \cite{Ichihara:2004} by way of asymptotic analysis on a dual representation of the nonlinear filtering equation. 
In these papers, the coefficients of the signal processes are assumed to be periodic. 
The approach in \cite{Ichihara:2004} is novel as the first application of backward stochastic differential equations for homogenization of Zakai-type stochastic partial differential equations (SPDEs). 

Convergence of the filter for a random ordinary differential equation with intermediate timescale and perturbed by a fast Markov process was investigated in \cite{Lucic:2003amo}. 
A two timescale problem with correlation between the slow process and observation process, but where the slow dispersion coefficient does not depend on the fast process, 
is investigated in \cite{Qiao:2019o}.
The main result is that the filter converges in $L^1$ sense to the lower dimensional filter. 
An energy method approach is used in \cite{Zhang:2019sd} to show that the probability density of the reduced nonlinear filtering problem approximates the original problem when the signal process has constant diffusion coefficients, periodic drift coefficients and the observation process is only dependent on the slow process. 

Convergence of the nonlinear filter is shown in a very general setting in \cite{Kleptsina:1997}, based on convergence in total variation distance of the law of $(\Xe, \Ye)$. 
In the examples of \cite{Kleptsina:1997}, the diffusion coefficient is not allowed to depend on the fast component.

The work of nonlinear filter approximation given in \cite[Chapter 6]{Kushner:1990}, is for a two timescale jump-diffusion process, but with no correlation between signal and observation process. 
The difference of the actual unnormalized conditional measure and the reduced conditional measure is shown to converge to zero in distribution.
Standard results then yield convergence in probability of the fixed time marginals. 
The method of proof is by averaging the coefficients of the SDEs for the unnormalized filters and showing that the limits of both filters satisfy the same SDE, which possess a unique solution.  
In \cite{Beeson:2020jma-arXiv} a similar approach to \cite[Chapter 6]{Kushner:1990} is used to study a broader multiple timescale correlation filtering problem where an intermediate scaling term exists and there is correlation between the slow and observation processes.
The authors make use of the perturbed test function approach where the correctors are solutions of Poisson equations to manage the difficulties introduced by the intermediate timescale. 
The main result in \cite{Beeson:2020jma-arXiv} is that $\pi^{\epsilon, x} \rightarrow \pi^0$ in probability for a metric generating the weak topology. 

In contrast to other papers on the convergence of the nonlinear filter for the multiple timescale problem, Imkeller et al. \cite{Imkeller:2013} showed a quantitative rate of convergence of $\epsilon$ for the system in Eq. \ref{equation: multiscale signal process}, but without intermediate timescale nor correlation of the slow process with the observation process. 
This is accomplished using a suitable asymptotic expansion of the dual of the Zakai equation and then harnessing a probabilistic representation of the SPDEs in terms of backward doubly stochastic differential equations. 
The approach of \cite{Imkeller:2013} is extended in this paper to cover the case of correlation between the observation process and the coarse-grain process. 
The analysis is therefore similar, with the exception of additional methods to handle the components of the dual of the Zakai equation due to the correlation and the final argument of the main proof. 

\begin{theorem*}[Main Result]
Under the assumptions stated in Theorem \ref{theorem: main result of quantitative convergence},
for every $p \geq 1$, $T \geq 0$, there exists a $C > 0$ such that for every $\varphi \in C^4_b(\R{m}; \R{})$,
\begin{align*}
\E[\Q{}]{ \left| \pi^{\epsilon, x}_T (\varphi) - \pi^0_T (\varphi) \right|^p }
\leq \epsilon^p C \norm{\varphi}{4, \infty}^p .
\end{align*}
In particular, there exists a metric $d$ on the space of probability measures on $\R{m}$, such that $d$ generates the topology of weak convergence, and such that for every $T \geq 0$, there exists $C > 0$ so that
\begin{align*}
\E[\Q{}]{ d( \pi^{\epsilon, x}_T , \pi^0_T ) }
\leq \epsilon C.
\end{align*}
\end{theorem*}

To prove the main result, we first setup the full problem in Section \ref{section: problem statement}, provide useful notation, and the main result with full details. 
The averaged SDE, Kushner-Stratonovich and Zakai equations are provided in this section as well. 
Having introduced the Zakai equations, we introduce their dual process representations in Section \ref{section: dual process} and explain how working with the dual process will allow us to prove the main result. 
In Section \ref{section: probabilistic representation of SPDEs}, a probabilistic representation of the dual processes is given.
Having established the necessary tools for the analysis, we provide preliminary estimates in Section \ref{section: preliminary estimates (chapter: quantitative)} and then the main analysis in Section \ref{section: main analysis of quantitative convergence}

\section{Problem Statement}
\label{section: problem statement}

In this section, we provide the full problem statement, some notation and the main result. 
We consider a filtered probability space $(\Omega, \mathcal{F}, (\mathcal{F}_t)_{t \geq 0}, \Q{})$ supporting a $(w + v + u)$-dimensional $\mF_t$-adapted Brownian motion $(W, V, U)$.
We will work with the following system of SDEs,
\EquationAligned{
\label{equation: original SDE problem setup}
dX^\epsilon_t 
&= b(X^\epsilon_t, Z^\epsilon_t) dt + \sigma(X^\epsilon_t, Z^\epsilon_t) dW_t, \\[1.5mm]
dZ^\eps_t 
&= \frac{1}{\eps^2} f(X^\eps_t, Z^\eps_t) dt + \frac{1}{\eps}g(X^\eps_t, Z^\eps_t) dV_t, \\[3mm]
dY^\eps_t 
&= 
h(X^\eps_t, Z^\eps_t) dt + \alpha dW_t + \gamma dU_t,  \quad \Ye_0 = 0 \in \R{d},\\
}
where $b : \R{m} \times \R{n} \rightarrow \R{m}$, $\sigma : \R{m} \times \R{n} \rightarrow \R{m} \times \R{w}$, $f : \R{m} \times \R{n} \rightarrow \R{n}$, $g : \R{m} \times \R{n} \rightarrow \R{n} \times \R{v}$ and $h : \R{m} \times \R{n} \rightarrow \R{d}$ are Borel measurable functions. 
The initial distribution of $(X, Z)$ is denoted by $\Q{}_{(\Xe_0, \Ze_0)}$ and is assumed independent of the $(W, V, U)$ Brownian motion.
$\Q{}_{(\Xe_0, \Ze_0)}$ is also assumed to have finite moments for all orders. 
In Eq. \ref{equation: original SDE problem setup}, $0 < \epsilon \ll 1$, is a timescale separation parameter.
We consider the case where $\alpha \in \R{d \by w}, \gamma \in \R{d \by u}$, and assume the following to be true
\begin{align*}
K \equiv \alpha \alpha^* + \gamma \gamma^* \succ 0, \quad \gamma \gamma^* \succ 0.
\end{align*}
This implies the existence of a unique $\R{d \by d} \ni \kappa \succ 0$ of lower triangular form, such that $K = \kappa \kappa^*$. 
Hence there exists a unique $\kappa^{-1}$, such that we can define an auxiliary observation process
\begin{align}
\label{equation: auxiliary observation process}
Y^{\epsilon, \kappa}_t = \int_0^t \kappa^{-1} d\Ye_s = \int_0^t \kappa^{-1} h(\Xe_s, \Ze_s) ds + B_t, \quad Y^{\epsilon, \kappa}_0 = 0 \in \R{d},
\end{align}
where
\begin{align*}
B_t = \kappa^{-1} \left( \alpha dW_t + \gamma dU_t \right), 
\end{align*}
is a standard $d$-dimensional Brownian motion under $\Q{}$. 

We are interested in the convergence of the $x$-marginal of the normalized filter, $\pi^{\epsilon, x}$, the conditional distribution of the signal given the observation filtration, to an averaged form. 
In particular, for any test function $\varphi \in C^2_b(\R{m} \times \R{n}; \R{})$ and time $t \in [0 , T]$, the normalized filter can be characterized as
\begin{align}
\label{equation: normalized conditional distribution as expectation}
\pi^\epsilon_t(\varphi) = \E[\Q{}]{ \varphi(\Xe_t, \Ze_t) }[\mYe_t], 
\end{align}
where $\mYe_t \equiv \sigma(\set{\Ye_s}[s \in [0, t]]) \vee \mathcal{N}$, the $\sigma$-algebra generated by the observation process over the interval $[0, t]$, joined with $\mathcal{N}$, the $\Q{}$ negligible sets. 

Because the filtrations generated by $\Ye$ and $Y^{\epsilon, \kappa}$ are equivalent, from the point of view of $\pi^\epsilon$ we can use either. 
Hence, let us redefine the sensor function $h \leftarrow \kappa^{-1} h$, the coefficients $\alpha \leftarrow \kappa^{-1} \alpha$ and $\gamma \leftarrow \kappa^{-1} \gamma$, so that the observation process can be redefined as
\begin{align}
\label{equation: redefined observation process}
dY^\eps_t = h(X^\eps_t, Z^\eps_t) dt + dB_t,  \quad \Ye_0 = 0 \in \R{d},
\end{align}
where $B = \alpha W + \gamma U$ is a standard Brownian motion under $\Q{}$ and still correlated with $W$. 

In Eq. \ref{equation: original SDE problem setup}, we identify the infinitesimal generators of the SDEs as follows, 
\begin{align*}
\mG_S (x, z) &\equiv \sum_{i=1}^m b_i(x, z) \frac{\partial}{\partial x_i} + \frac{1}{2} \sum_{i, j=1}^m (\sigma \sigma^*)_{ij} (x, z) \frac{\partial^2}{\partial x_i \partial x_j}, \\
\mG_F (x, z) &\equiv \sum_{i=1}^n f_i(x, z) \frac{\partial}{\partial z_i} + \frac{1}{2} \sum_{i, j=1}^n (gg^*)_{ij}(x, z) \frac{\partial^2}{\partial z_i \partial z_j},  \\
\mG^\epsilon &\equiv \frac{1}{\epsilon^2} \mG_F + \mG_S.
\end{align*}
The Kushner-Stratonovich equation for the time evolution of the filter $\pi^\epsilon$, acting on a test function $\varphi \in C^2_b(\R{m} \by \R{n}; \R{})$, is
\EquationAligned{\label{equation: kushner-stratonovich}
\pi^\epsilon_t(\varphi) &= \pi^\epsilon_0(\varphi) + \int_0^t \pi^\epsilon_s (\mG^\epsilon \varphi) ds + \int_0^t \langle \pi^\epsilon_s ( \varphi h + \alpha \sigma^* \nabla_x \varphi ) - \pi^\epsilon_s (\varphi) \pi^\epsilon_s (h), d\Ye_s - \pi^\epsilon_s (h) ds \rangle, \\
\pi^\epsilon_0(\varphi) &= \E[\Q{}]{ \varphi(\Xe_0, \Ze_0) }.
}
When we are interested in estimating test functions of $\Xe$ only, i.e., $\varphi \in C^2_b(\R{m}; \R{})$, we consider the $x$-marginal of $\pi^\epsilon$, 
\begin{align}
\label{equation: marginal conditional distribution as expectation}
\pi^{\epsilon, x}_t(\varphi) = \int \varphi(x) \pi^\epsilon_t(dx, dz).
\end{align}

\subsection{Diffusion Approximation and the Averaged Filter}
\label{subsection: diffusion approximation and averaged filtering equations (chapter: quantitative)}

The theory of homogenization of stochastic differential equations shows that if the process $Z^{\epsilon, x}$, 
\begin{align}
\label{equation: SDE, fast process with fixed slow state}
dZ^{\eps, x}_t = \frac{1}{\eps^2} f(x, Z^{\eps, x}_t) dt + \frac{1}{\eps}g(x, Z^{\eps, x}_t) dV_t,
\end{align}
 is ergodic with stationary distribution $\mu_\infty(x)$, then under appropriate conditions, in the limit $\epsilon \rightarrow 0$ the process $\Xe$ converges in distribution to a Markov process $X^0$ with infinitesimal generator 
\begin{align}
\label{equation: generator for the averaged SDE process}
\overbar{\mG_S} (x) &\equiv \sum_{i = 1}^m \overbar{b}_i(x) \frac{\partial}{\partial x_i} + \frac{1}{2} \sum_{i, j = 1}^m \overbar{a}_{ij}(x) \frac{\partial^2}{\partial x_i \partial x_j},
\end{align}
where the averaged drift and diffusion coefficients are
\begin{align*}
\overbar{b}(x) \equiv \int_{\R{n}} b(x, z) \mu_\infty(dz; x), 
\qquad \text{and} \qquad
\overbar{a}(x) \equiv \int_{\R{n}} a(x, z) \mu_\infty(dz; x).
\end{align*}
Here we denote the diffusion coefficient $a = \sigma \sigma^*$.
Additionally, let us define 
\begin{align*}
\overbar{h}(x) \equiv \int_{\R{n}} h(x, z) \mu_\infty(dz; x), 
\qquad \text{and} \qquad
\overbar{\sigma}(x) \equiv \int_{\R{n}} \sigma(x, z) \mu_\infty(dz; x).
\end{align*}

The aim of this paper is to show that the $x$-marginal filter $\pi^{\epsilon, x}$ can be approximated by an averaged filter $\pi^0$. 
We will show the existence and uniqueness of $\pi^0$ in Section \ref{section: the averaged conditional distributions}. 
This is done by defining $\pi^0$ from the Kallianpur-Striebel formula and the existence and uniqueness of an unnormalized averaged filter $\rho^0$. 
The averaged filter will depend on the averaged coefficients $\overbar{b}, \overbar{a}, \overbar{\sigma}$ and $\overbar{h}$.

\subsection{Notation and Main Theorem}

Before stating the main result of the paper, we set a few definitions and assumptions that will be used throughout the paper. 
We will use $\N{}_0$ to denote $\{ 0, 1, 2, \hdots \}$ and $\N{}$ for $\set{1, 2, \hdots}$.
Let $H_f$ denote the assumption that there exists a constant $C > 0$, exponent $\alpha > 0$ and an $R > 0$ such that for all $|z| > R$,  
\begin{align}
\tag{$H_f$}
\label{assumption: positive recurrence}
\sup_{x \in \R{m}} \langle f(x, z), z \rangle \leq - C |z|^\alpha.
\end{align}
\ref{assumption: positive recurrence} is a recurrence condition, which provides the existence of a stationary distribution, $\mu_\infty(x)$, for the process $Z^x$.
Let $H_g$ denote the assumption that there are $0 < \lambda \leq  \Lambda < \infty$, such that for any $(x, z) \in \R{m} \times \R{n}$,
\begin{align}
\tag{$H_g$}
\label{assumption: uniform ellipticity}
\lambda I \preceq g g^*(x, z) \preceq \Lambda I, 
\end{align}
where $\preceq$ is the order relation in the sense of positive semidefinite matrices. 
\ref{assumption: uniform ellipticity} is a uniform ellipticity condition, which provides the uniqueness of the stationary distribution. 
We will say that a function $\theta: \R{m} \times \R{n} \rightarrow \R{}$ is centered with respect to $\mu_\infty(x)$, if for each $x$ 
\[
\int \theta(x, z) \mu_\infty(dz; x) = 0, \quad \forall x \in \R{m}. 
\]
If $\varphi(x, z) \in C^{k, l}_b (\R{m} \times \R{n}; \R{n})$, then $\varphi$ is $k$-times continuously differentiable in the $x$-component, $l$-times continuously differentiable in the $z$-component, and all partial derivatives $\partial^{l'}_z \partial^{k'}_x \varphi$ for $0 \leq k' \leq k$, $0 \leq l' \leq l$ are bounded.
Let $HF^{k, l}$ for $k, l \in \N{}_0$ denote the following assumption:
\hypertarget{HF}{}
\begin{align}
\tag{$HF^{k, l}$}
\label{assumption: regularity and boundedness of fast coefficients}
f \in C^{k, l}_b(\R{m} \times \R{n}; \R{n}) 
\qquad
\text{and}
\qquad
g \in C^{k, l}_b(\R{m} \times \R{n}; \R{n \times k}).
\end{align}
Similarly, let $HS^{k, l}$ for $k, l \in \N{}_0$ denote the assumption:
\hypertarget{HS}{}
\begin{align}
\tag{$HS^{k, l}$}
\label{assumption: regularity and boundedness of slow coefficients}
b \in C^{k, l}_b(\R{m} \times \R{n}; \R{m})
\qquad
\text{and}
\qquad
\sigma \in C^{k, l}_b(\R{m} \times \R{n}; \R{m \times k}),
\end{align}
and $HO^{k, l}$ for $k, l \in \N{}_0$ denote the assumption:
\hypertarget{HO}{}
\begin{align}
\tag{$HO^{k, l}$}
\label{assumption: regularity and boundedness of observation function}
h \in C^{k, l}_b(\R{m} \times \R{n}; \R{d}).
\end{align}
We use the notation $k = (k_1, \hdots, k_m) \in \mathbb{N}^m_0$ for a multiindex with order $|k| = k_1 + \hdots + k_m$ and define the differential operator
\begin{align*}
D^k_x = \frac{\partial^{|k|}}{\partial {x_1}^{k_1} \hdots \partial x_m^{k_m}}.
\end{align*}
Lastly, the relation $a \lesssim b$ will indicate that $a \leq C b$ for a constant $C > 0$ that is independent of $a$ and $b$, but that may depend on parameters that are not critical for the bound being computed.

Having introduced the necessary definitions and equations, we now state the main result fully.
\begin{theorem}
\label{theorem: main result of quantitative convergence}
Assume \ref{assumption: positive recurrence}, \ref{assumption: uniform ellipticity}, \hyperlink{HF}{$HF^{8, 4}$}, $b \in C^{7, 4}_b$, $\sigma \in C^{8, 4}_b$, and \hyperlink{HO}{$HO^{8, 4}$}.
Additionally, assume that the initial distribution $\Q{}_{(\Xe_0, \Ze_0)}$ has finite moments of every order. 
Then for any $p \geq 1, T \geq 0$ we have that for every $\varphi \in C^4_b(\R{m}; \R{})$,
\begin{align*}
\E[\Q{}]{ \left| \pi^{\epsilon, x}_T (\varphi) - \pi^0_T (\varphi) \right|^p }
\lesssim \epsilon^p \norm{\varphi}{4, \infty}^p .
\end{align*}
Further, there exists a metric $d$ on the space of probability measures on $\R{m}$ that generates the topology of weak convergence, such that 
\begin{align*}
\E[\Q{}]{ d( \pi^{\epsilon, x}_T , \pi^0_T ) }
\lesssim \epsilon .
\end{align*}

\begin{proof}
The proof of the first result is given by Corollary \ref{corollary: moment estimate of normalized measure error}.
The proof of the second result is from Lemma \ref{lemma: convergence of normalized measure error under metric generating weak topology}.
\end{proof}
\end{theorem}

Before moving beyond this theorem statement, we provide some quick remarks. 

{
\rem{
From $\lim_{\epsilon \rightarrow 0 } \E[\Q{}] { d( \pi^{\epsilon, x}_T , \pi^0_T ) } = 0$, we retrieve convergence in probability,
\begin{align*}
\lim_{\epsilon \rightarrow 0} \Q{} \left( d( \pi^{\epsilon, x}_T , \pi^0_T ) \geq \delta \right) 
\leq \frac{1}{\delta} \lim_{\epsilon \rightarrow 0 } \E[\Q{}] { d( \pi^{\epsilon, x}_T , \pi^0_T ) } 
= 0, \qquad \text{for each $\delta > 0$}. 
\end{align*}
And by the Borel-Cantelli lemma we can choose $(\epsilon_n)$ so that $\pi^{\epsilon, x}_n$ will a.s. converge weakly to $\pi^0$. 
}
}

{
\rem{
Some quick comparisons to the main result in \cite{Imkeller:2013}.
There the scaling for the fast process was of order one, whereas in this paper we use order two. 
Therefore, the rate of convergence is the same in the two works. 
The only difference in the conditions of our Theorem \ref{theorem: main result of quantitative convergence} and the equivalent one in \cite{Imkeller:2013}, is that we require $\sigma \in C^{8, 4}_b$ instead of $C^{7, 4}_b$. 
This extra regularity in the slow component of the function is due to the correlation between the slow process and observation process, which then appears in our backward stochastic differential equations of Section \ref{section: main analysis of quantitative convergence}.
}
}

\subsection{Change of Probability Measure and the Zakai Equation}

In Section \ref{section: main analysis of quantitative convergence}, we will be interested in working with the unnormalized conditional measure since it will satisfy a linear evolution equation. 
To define the unnormalized conditional measure requires a change of probability measure transformation, which we will perform for each $\epsilon$.  
Let us denote the new collection of probability measures by $(\Pe)$. 
For any fixed $\epsilon$, $\Pe$ and $\Q{}$ will be mutually absolutely continuous with Radon-Nikodym derivatives 
\begin{align*}
D^\eps_t 
\equiv \restrict{\frac{d\P^\eps}{d\Q{}}}{\mF_t} 
= \exp \left( - \int_0^t \langle h(X^\eps_s, Z^\eps_s), dB_s \rangle - \frac{1}{2} \int_0^t \left| h(X^\eps_s, Z^\eps_s) \right|^2 ds \right),
\end{align*}
\begin{align*}
\wt{D}^\eps_t 
\equiv (D^\eps_t)^{-1} 
= \restrict{\frac{d\Q{}}{d\P^\eps}}{\mF_t} 
= \exp \left( \int_0^t \langle h (\Xe_s, \Ze_s), d\Ye_s \rangle - \frac{1}{2} \int_0^t \left| h (X^\eps_s, Z^\eps_s) \right|^2 ds \right). 
\end{align*}
Then by Girsanov's theorem, under $\P^\eps$ the process $\Ye$ is a Brownian motion.
For a fixed test function $\varphi \in C^2_b(\R{m} \times \R{n}; \R{})$ and time $t \in [0 , T]$, we characterize the unnormalized conditional measure $\rho^\epsilon_t$ as, 
\begin{align*}
\rho^\eps_t(\varphi) 
= \E[\Pe]{\varphi(\Xe_t, \Ze_t)\wt{D}^\eps_t \st \mYe_t},
\end{align*}
and its relation to $\pi^\epsilon$ is given by the Kallianpur-Striebel formula, 
\begin{align*}\label{equation: Kallianpur-Striebel}
\pi^\eps_t(\varphi) 
= \frac{\E[\Pe]{\varphi(\Xe_t, \Ze_t)\wt{D}^\eps_t \st \mYe_t}}{\E[\Pe]{ \wt{D}^\eps_t \st \mYe_t}} 
= \frac{\rho^\eps_t(\varphi)}{\rho^\eps_t(1)},
\qquad \forall t \in [0, T], 
\qquad \text{$\Q{}, \Pe$-a.s}.
\end{align*}
The action of $\rho^\epsilon$ on test functions $\varphi \in C^2_b(\R{m} \by \R{n}; \R{})$ gives the Zakai evolution equation, 
\EquationAligned{
\label{equation: Zakai}
\rho^\epsilon_t(\varphi) &= \rho^\epsilon_0(\varphi) + \int_0^t \rho^\epsilon_s \left( \mG^\epsilon \varphi \right) ds + \int_0^t \langle \rho^\epsilon_s ( \varphi h + \alpha \sigma^* \nabla_x \varphi ),  d\Ye_s \rangle, \\
\rho^\epsilon_0(\varphi) &= \E[\Q{}]{ \varphi(\Xe_0,\Ze_0) }.
}
When $\varphi \in C^2_b(\R{m}; \R{})$, we consider the $x$-marginal, 
\begin{align*}
\rho^{\epsilon, x}_t(\varphi) = \int \varphi(x) \rho^\epsilon_t(dx, dz),
\end{align*}
which is related to $\pi^{\epsilon, x}$ through the Kallianpur-Striebel formula, 
\begin{align*}
\pi^{\epsilon, x}_t(\varphi) 
= \frac{\rho^{\epsilon, x}_t(\varphi)}{\rho^{\epsilon, x}_t(1)},
\qquad \forall t \in [0, T], 
\qquad \text{$\Q{}, \Pe$-a.s}.
\end{align*}

\subsubsection{The Averaged Conditional Distributions}
\label{section: the averaged conditional distributions}

In this section, we show that there exists a probability measure-valued process $\pi^0$, the averaged filter, that is defined in terms of a measure-valued process $\rho^0$, the averaged unnormalized filter, from the Kallianpur-Striebel formula. 
To do so, we start by defining under $\Q{}$, the SDE $\XOe$ satisfying the equation,  
\begin{align*}
d\XOe_t 
&= \overbar{b}(\XOe_t) dt + ( \overbar{a}(\XOe_t) - \overbar{\sigma} \overbar{\sigma}^*(\XOe_t) )^{1/2} d\widehat{W}_t 
+ \overbar{\sigma}(\XOe_t) \left( dW_t + \alpha^* h(\Xe_t, \Ze_t) dt - \alpha^* \overbar{h}(\XOe_t) dt \right), \\
\XOe_0 
&\sim \Q{}_{\Xe_0},
\end{align*}
which under the change of measure to $\Pe$ becomes
\EquationAligned{
\label{equation: auxiliary averaged SDE under the change of measure}
d\XOe_t 
&= \overbar{b}(\XOe_t) dt + ( \overbar{a}(\XOe_t) - \overbar{\sigma} \overbar{\sigma}^*(\XOe_t) )^{1/2} d\widehat{W}_t 
+ \overbar{\sigma}(\XOe_t) \left( d\overbar{W}_t - \alpha^* \overbar{h}(\XOe_t) dt \right), \\
\XOe_0 
&\sim \Q{}_{\Xe_0},
}
where $\overbar{W}_t$ is a standard Brownian motion under $\Pe$. 
The Cholesky factor $( \overbar{a}(X^0_t) - \overbar{\sigma} \overbar{\sigma}^*(X^0_t) )^{1/2}$ exists, since from an application of Jensen's inequality $\overbar{a}(x) - \overbar{\sigma} \overbar{\sigma}^*(x) \succeq 0$ for each $x \in \R{m}$. 
We now define under $\Pe$, the process
\begin{align*}
\wt{D}^0_t 
= \exp \left( \int_0^t \langle \overbar{h} (\XOe_s), d\Ye_s \rangle - \frac{1}{2} \int_0^t \left| \overbar{h} (\XOe_s) \right|^2 ds \right),
\end{align*}
which also satisfies the relation,
\begin{align*}
\wt{D}^0_t 
= 1 + \int_0^t \langle \wt{D}^0_s \overbar{h} (\XOe_s), d\Ye_s \rangle.
\end{align*}

\begin{lemma}
\label{lemma: existence of averaged unnormalized filter}
Assume $\overbar{h}$ is a bounded function. 
Then there exists a measure-valued process $(\rho^0_t)_{t \geq 0}$ such that for all $\varphi \in C_b$,
\begin{align*}
\rho^0_t(\varphi) = \E[\Pe]{\varphi(\XOe_t) \wt{D}^0_t}[\mYe_t], \qquad \Pe\text{-a.s.}
\end{align*}
Additionally, for every $\varphi \in C^2_b(\R{m}; \R{})$, $\rho^0(\varphi)$ satisfies the equation,
\EquationAligned{
\label{equation: homogenized Zakai}
\rho^0_t(\varphi) &= \rho^0_0(\varphi) + \int_0^t \rho^0_s(\overbar{\mG_S} \varphi) ds + \int_0^t \langle \rho^0_s(\varphi \overbar{h} + \alpha \overbar{\sigma}^* \nabla_x \varphi ),  d\Ye_s \rangle, \\
\rho^0_0(\varphi) &= 
\int \varphi(x) \Q{}_{\Xe_0}(dx),
}
where $\Q{}_{\Xe_0}$ is the initial distribution of $\Xe$ and the solution is unique if the coefficients of $\overbar{b}, \overbar{a}, \overbar{\sigma}$, and $\overbar{h}$ are in $C^3_b$.
\begin{proof}
Because $\overbar{h}$ is bounded, we have by the same proof as Lemma \ref{lemma: moments of Girsanov (chapter: quantitative)}, the uniform bound 
\[
\sup_{\epsilon \in (0, 1]} \sup_{t \leq T} \E*[\Pe]{ |\wt{D}^0_t|^p } < \infty,
\]
for $p \geq 2, T > 0$.
Then $\E*[\Pe]{ \wt{D}^0_t } = 1$ and therefore $\Qo(\cdot) = \int_{\cdot} \wt{D}^0_t(\omega) \Pe(d\omega)$ is a new probability measure. 
Since $\wt{D}^0_t$ is $\Pe$-a.s. strictly positive, we also have that $\Qo$ is equivalent to $\Pe$. 
By the Kallianpur-Striebel formula we know that the following holds, 
\begin{align}
\label{equation: KS for averaged unnormalized filter}
\frac{\E[\Pe]{\varphi(\XOe_t) \wt{D}^0_t}[\mYe_t]}{\E[\Pe]{\wt{D}^0_t}[\mYe_t]} 
= \E[\Qo]{\varphi(\XOe_t)}[\mYe_t].
\end{align}
Indeed, for any $\mYe_t$-measurable random variable $\xi$: 
\begin{align*}
\E[\Pe]{\xi \E[\Qo]{\varphi(\XOe_t)}[\mYe_t] \E[\Pe]{\wt{D}^0_t}[\mYe_t]}
&= \E[\Pe]{ \xi \E[\Qo]{\varphi(\XOe_t)}[\mYe_t] \wt{D}^0_t } \\
&= \E[\Qo]{ \xi \E[\Qo]{\varphi(\XOe_t)}[\mYe_t] } \\
&= \E[\Qo]{ \xi \varphi(\XOe_t) } 
= \E[\Pe]{ \xi \varphi(\XOe_t) \wt{D}^0_t},
\end{align*}
and therefore 
\begin{align*}
\E[\Qo]{\varphi(\XOe_t)}[\mYe_t] \E[\Pe]{\wt{D}^0_t}[\mYe_t] = \E[\Pe]{ \varphi(\XOe_t) \wt{D}^0_t}[\mYe_t]. 
\end{align*}
As a consequence of $\wt{D}^0_t > 0$ $\Pe$-a.s., we have that the random variable $\E[\Pe]{\wt{D}^0_t}[\mYe_t]$ is $\Pe$-a.s. strictly positive, so we can divide by the variable to obtain Eq. \ref{equation: KS for averaged unnormalized filter}. 
Therefore, there exists a regular $\Qo$-conditional probability $\pi^0_t$ such that $\Qo$-a.s. $\pi^0_t(\varphi) = \E[\Qo]{\varphi(\XOe_t)}[\mYe_t]$ for each fixed $t$. 
Since $\Pe$ and $\Qo$ are equivalent, the identity also holds $\Pe$-a.s.
This only gives a random measure at one fixed time, but \cite[Theorem 2.24, p.29]{Bain:2009wo} can be used to obtain $\pi^0$ as a probability measure-valued process.
Note that in our setup, $\XOe$ is not a Markov process, as is assumed in \cite[Theorem 2.24, p.29]{Bain:2009wo}, but this is not important for that theorem and one could always consider $\XOe$ to be a component of the Markov process $(\XOe, \Xe, \Ze)$. 
We now set $\rho^0_t(\varphi) = \pi^0_t(\varphi) \E[\Pe]{ \wt{D}^0_t }[\mYe_t]$ to get the first part of the proof. 
For the last part, by standard construction of the Zakai equation (see for instance \cite{Bain:2009wo}), $\rho^0$ satisfies Eq. \ref{equation: homogenized Zakai} and uniqueness follows from \cite[Theorem 3.1, p.454]{Rozovskii:1991}.
\end{proof}
\end{lemma}

Therefore based on Lemma \ref{lemma: existence of averaged unnormalized filter}, the averaged (normalized) filter $\pi^0$ is then related to $\rho^0$ by the Kallianpur-Striebel relation, 
\begin{align}
\label{equation: homogenized filter via Kallianpur-Streibel}
\pi^0_t(\varphi) = \frac{\rho^0_t(\varphi) }{\rho^0_t(1) },
\qquad \forall t \in [0, T], \quad \forall \varphi \in C_b(\R{m}; \R{}).
\end{align}

{
\rem{
An interesting observation regarding Eq. \ref{equation: homogenized Zakai}, is that we may have $\overbar{\sigma} = 0$, and this implies that the SDE for the averaged filter may have no correlation. 
}
}

{
\rem{
Note that $\pi^0$ is not the filter for the averaged system, and hence $\rho^0$ is also not the unnormalized conditional measure for the averaged system, which would instead satisfy the following equation, 
\EquationAligned*{
\label{equation: Zakai for averaged system}
\overbar{\rho}_t(\varphi) &= \overbar{\rho}_0(\varphi) + \int_0^t \overbar{\rho}_s(\overbar{\mG_S} \varphi) ds + \int_0^t \langle \overbar{\rho}_s(\varphi \overbar{h} + \alpha \sqrt{\overbar{a}}^* \nabla_x \varphi ),  d\overbar{Y}_s \rangle, \\
\overbar{\rho}_0(\varphi) &= 
\int \varphi(x) \Q{}_{\Xe_0}(dx),
}
where 
\begin{align*}
\overbar{Y}_t = \int_0^t \overbar{h}(X^0_s) ds + B_t,
\end{align*}
and $X^0$ is a diffusion process with infinitesimal generator $\overbar{\mG_S}$ under $\Q{}$.
}
}

\section{Dual Process to the Unnormalized Conditional Distribution}
\label{section: dual process}

We now introduce an idea by \cite{Pardoux:1980a} that is an important transition to the method of proof used in this paper. 
The idea is to define a function-valued process for any fixed $\varphi \in C^2_b(\R{m}; \R{})$.
The function-valued process will be the dual of $\rho^\epsilon$ in an appropriate sense. 
We first define
\begin{align*}
\wDe_{t, T} 
= \exp \left( \int_t^T \langle h (\Xe_s, \Ze_s), d\Ye_s \rangle - \frac{1}{2} \int_t^T \left| h (X^\eps_s, Z^\eps_s) \right|^2 ds \right),
\end{align*}
which is $\Q{}$-a.s. equal to $\wDe_T  (\wDe_t)^{-1}$. 
Fixing $\varphi \in C^2_b(\R{m}; \R{})$, we then define the dual process at time $t \in [0, T]$ as
\begin{align*}
v^{\eps, T, \varphi}_t(x, z)
\equiv \E[\Pe_{t, x, z}]{ \varphi(\Xe_T) \wDe_{t, T} }[\mYe_{t, T} ],
\end{align*}
where $\Pe_{t, x, z}$ is the change of probability measure that results when $(\Xe_s, \Ze_s)$ takes the constant value $(x, z)$ for $s \in [0, t]$ and then follows the dynamics given by Eq. \ref{equation: multiscale signal process} for $s > t$. 
$v^{\eps, T, \varphi}_t(x, z)$ is called the dual process, because for any $t \in [0, T]$ we have 
\begin{align*}
\rho^{\epsilon, x}_T(\varphi) 
= \rho^\epsilon_t ( v^{\eps, T, \varphi}_t ), \qquad \Pe\text{-a.s.} 
\end{align*}
This also means that $\rho^{\epsilon, x}_T(\varphi)  =  \rho^\epsilon_0 ( v^{\eps, T, \varphi}_0 )$ and therefore
\begin{align*}
\rho^{\epsilon, x}_T(\varphi)  
= \int_{\R{m} \times \R{n}} v^{\eps, T, \varphi}_0(x, z) \Q{}_{(\Xe_0, \Ze_0)}(dx, dz).
\end{align*}

We can similarly define the dual process for $\rho^0$.
Following the construction in Lemma \ref{lemma: existence of averaged unnormalized filter}, we would have 
\begin{align*}
v^{0, T, \varphi}_t(x) \equiv \E[\Pe_{t, x}]{\varphi(\XOe_T) \wt{D}^0_{t, T} \st \mYe_{t, T}},
\end{align*}
with the same property that $\rho^0_T(\varphi) = \rho^0_0 ( v^{0, T, \varphi}_0 )$. 
Again $\Pe_{t, x}$ is the change of probability measure that results from $\XOe_s$ taking the constant value $x$ for $s \in [0, t]$ and then follows the dynamics given by the SDE in Eq. \ref{equation: auxiliary averaged SDE under the change of measure}.
The definition of $\wt{D}^0_{t, T}$ in $v^{0, T, \varphi}_t$ is 
\begin{align*}
\wt{D}^0_{t, T} 
= \exp \left( \int_t^T \langle \overbar{h} (\XOe_s), d\Ye_s \rangle - \frac{1}{2} \int_t^T \left| \overbar{h} (\XOe_s) \right|^2 ds \right),
\end{align*}
which is $\Pe{}$-a.s. equal to $\wt{D}^0_{t, T}  = \wt{D}^0_T ( \wt{D}^0_t )^{-1} $.

\subsection{The Dual Process and Filter Convergence}
\label{subsection: the dual process and filter convergence}

We now show the usefulness of the dual process in showing the convergence of $\rho^{\epsilon, x} \rightarrow \rho^0$. 
We again fix $\varphi \in C^2_b(\R{m}; \R{})$ and $p \geq 1$. 
Then from Jensen's inequality and Fubini's theorem we have the following relation,
\EquationAligned{
\label{equation: inequality, dual residual implies zakai residual}
\E[\Pe]{ \left| \rho^{\epsilon, x}_T(\varphi) - \rho^0_T(\varphi) \right|^p } 
=  \E[\Pe]{ \left| \int v^{\epsilon, T, \varphi}_0(x, z) - v^{0, T, \varphi}_0(x) \Q{}_{(\Xe_0, \Ze_0)}(dx, dz) \right|^p }  \\
\leq  \E[\Pe]{ \int \left| v^{\epsilon, T, \varphi}_0(x, z) - v^{0, T, \varphi}_0(x) \right|^p \Q{}_{(\Xe_0, \Ze_0)}(dx, dz) }  \\
=  \int \E[\Pe]{ \left| v^{\epsilon, T, \varphi}_0(x, z) - v^{0, T, \varphi}_0(x) \right|^p } \Q{}_{(\Xe_0, \Ze_0)}(dx, dz).
}
This implies that if $\Q{}_{(\Xe_0, \Ze_0)}(dx, dz)$ is well behaved (e.g., finite moments of every order) then convergence of the $p$-th moment of $v^{\epsilon, T, \varphi}_0(x, z) - v^{0, T, \varphi}_0(x)$ to zero will imply convergence of  the $p$-th moment of $\rho^{\epsilon, x}_T(\varphi) - \rho^0_T(\varphi)$ to zero. 
Without loss of generality, we assumed that $\XOe$ had the same initial distribution as $\Xe$ in Section \ref{section: the averaged conditional distributions}, and hence why the integration is against $\Q{}_{(\Xe_0, \Ze_0)}$ in Eq. \ref{equation: inequality, dual residual implies zakai residual}.

\subsection{Evolution Equations for the Dual Process}

To introduce the next step in the techniques to prove convergence of the marginalized filter to the reduced order filter, we need to state the evolution equations for the dual processes $v^{\epsilon, T, \varphi}$ and $v^{0, T, \varphi}$. 
Both processes satisfy backward stochastic partial differential equations (BSPDE).
To facilitate the reading, we use $v^\epsilon$ and $v^0$ instead of the more verbose $v^{\epsilon, T, \varphi}$ and $v^{0, T, \varphi}$ in most of what follows. 
When clarity is needed, we will use the explicit notation.

The evolution equation for $v^\epsilon$ is given by  
\begin{align}
\label{equation: BSPDE, dual process}
 -dv^\eps_t 
= \mG^\eps v^\epsilon_t dt + \langle v^\epsilon_t h + \alpha \sigma^* \nabla_x v^\epsilon_t , d \overleftarrow{B}_t \rangle, \qquad v^\eps_T = \varphi,
\end{align}
where $\ds \int d\overleftarrow{B}_t$ will denote the backward It\^{o} integral.
The process $v^0$ is given by 
\begin{align}
\label{equation: BSPDE, averaged dual process}
 -dv^0_t 
= \overbar{\mG_S} v^0_t dt + \langle v^0_t \overbar{h} + \alpha \overbar{\sigma}^* \nabla_x v^0_t , d \overleftarrow{B}_t \rangle, \qquad v^0_T = \varphi.
\end{align}

\subsection{Expansion of the Dual Process}

Because $v^\epsilon$ satisfies a linear equation, we consider an expansion of $v^\epsilon$ using $v^0$ and a corrector $\psi$ and remainder $R$ term,
\begin{align*}
v^\epsilon_t(x, z) = v^0_t(x) + \psi_t(x, z) + R_t(x, z).
\end{align*}
Using this expansion in Eq. \ref{equation: BSPDE, dual process} and introducing terms for Eq. \ref{equation: BSPDE, averaged dual process}, we define $\psi$ and $R$ to satisfy the following linear BSPDEs
\begin{align}
-d \psi_t 
&= \left[ \frac{1}{\epsilon^2} \mG_F \psi_t + (\mG_S - \overbar{\mG_S}) v^0_t \right] dt + \langle v^0_t (h - \overbar{h}) , d \B \rangle + \langle \alpha (\sigma - \overbar{\sigma})^* \nabla_x v^0_t , d \B \rangle, 
&\qquad \psi_T = 0, \nonumber \\
-dR_t 
&= (\mG^\epsilon R_t + \mG_S \psi_t) dt + \langle (\psi_t +  R_t) h , d \B \rangle + \langle \alpha \sigma^* \nabla_x (\psi_t + R_t) , d \B \rangle,
&\qquad R_T = 0. \label{equation: BSPDE, remainder term}
\end{align}

Therefore to show convergence of the difference $v^\epsilon - v^0$, we can equivalently show convergence of $\psi$ and $R$ to zero as $\epsilon \rightarrow 0$: 
\begin{align*}
\E[\Pe]{ \left| v^{\epsilon, T, \varphi}_t(x, z) - v^{0, T, \varphi}_t(x) \right|^p }
= \E[\Pe]{ \left| \psi_t(x, z) + R_t(x, z) \right|^p } 
\lesssim \E[\Pe]{ \left| \psi_t(x, z) \right|^p } + \E[\Pe]{ \left| R_t(x, z) \right|^p }. 
\end{align*}
This will be our strategy in Section \ref{section: main analysis of quantitative convergence}.

\section{Probabilistic Representation of Stochastic PDEs}
\label{section: probabilistic representation of SPDEs}

We will now show that we can find a probabilistic representation of the dual processes. 
This representation will be given by backward doubly stochastic differential equations (BDSDEs), which are a generalization of the Feynman-Kac solution for semilinear second order parabolic SPDEs. 
First let us state a result for the classical solution of the dual processes, which are linear second order parabolic SPDEs of the general form: 
\begin{align}
\label{equation: general form of BSPDE}
-d \psi(\omega, t, x) 
= \mathcal{L} \psi (\omega, t, x) dt 
+ f (\omega, t, x) dt 
&+ \langle g (\omega, t, x) + G (\omega, t, x) \psi  (\omega, t, x) , d \B \rangle \\
&+ \langle F (\omega, t, x) \nabla_x \psi  (\omega, t, x) , d \B \rangle, \nonumber \\
\psi(T, x)
&= \varphi(\omega, x), \nonumber
\end{align}
where $\psi : \Omega \times [0, T] \times \R{m} \rightarrow \R{}$, $f: \Omega \times [0, T] \times \R{m} \rightarrow \R{}$, $g: \Omega \times [0, T] \times \R{m} \rightarrow \R{d}$, $G: \Omega \times [0, T] \times \R{m} \rightarrow \R{d}$, $F: \Omega \times [0, T] \times \R{m} \rightarrow \R{d \times m}$ and $\varphi: \Omega \times \R{m} \rightarrow \R{}$ are all jointly measurable, and $\B$ is a $d$-dimensional standard backward Brownian motion.
The generator given in Eq. \ref{equation: general form of BSPDE} has the form
\begin{align*}
\mathcal{L}(x) = \sum_{i = 1}^m b_i(x) \frac{\partial}{\partial x_i} + \frac{1}{2} \sum_{i, j = 1}^m a_{ij}(x) \frac{\partial^2}{\partial x_i \partial x_j},
\end{align*}
where $b: \R{m} \rightarrow \R{}$ and $a: \R{m} \rightarrow \mathbb{S}^{m \times m}$ are measurable ($\mathbb{S}^{m \times m}$ denotes the space of symmetric positive semidefinite matrices). 

Once we have stated the results on BSPDEs, we will then give the BDSDE representation in Section \ref{subsection: BDSDE}.
We will need some definitions for the necessary conditions on the classical solution of the SPDEs. 
Let us state those now, starting with the definition for the filtration $\mF^{0, B}_{t, s}$: let $0 \leq t \leq s \leq T$, 
\begin{align*}
\mF^{0, B}_{t, s}
= \sigma( \set{ B_u - B_t }[ t \leq u \leq s ] ) , 
\end{align*}
and let $\mF^B_{t, s}$ be the completion of $\mF^{0, B}_{t, s}$ under $\Pe$. 
We next define the space of adapted random fields of polynomial growth \hypertarget{P}{$\mathcal{P}_T(\R{m}; \R{n})$}: 
{
\defn{
\label{definition: random fields of polynomial growth}
$\mathcal{P}_T(\R{m}; \R{n})$ is the space of random fields of polynomial growth
\begin{align*}
H : \Omega \times [0, T] \times \R{m} \rightarrow \R{n}
\end{align*}
that are jointly measurable in $(\omega, t, x)$ and for fixed $(t, x)$, $\omega \mapsto H(\omega, t, x)$ is $\F$-measurable. 
Further, for fixed $\omega$ outside a null set, $H$ has to be jointly continuous in $(t, x)$, and it has to satisfy the following inequality: 
For every $p \geq 1$ there is $C_p, q > 0$, such that for all $x \in \R{m}$, 
\begin{align*}
\E{ \sup_{0 \leq t \leq T} | H(t, x) |^p } \leq C_p (1 + |x|^q ). 
\end{align*}
}
}

We denote with \hypertarget{DK}{$D^k$} a definition concerning conditions on the coefficients of the generator $\mathcal{L}$ of the BSPDE:
{
\defn{
\label{definition: conditions on coefficients of BSPDE generator}
We define the condition $D^k$ to indicate that $b \in C^k_b(\R{m}; \R{n})$, $a \in C^k_b(\R{m}; \mathbb{S}^{m \times m})$, and $a$ is degenerate elliptic: 
For every $\xi \in \R{m}$ and every $x \in \R{m}$, 
\begin{align*}
\langle a(x) \xi , \xi \rangle = \sum_{i, j = 1}^m a_{ij}(x) \xi_i \xi_j \geq 0,
\end{align*}
or succintly $a \succeq 0$.
}
}

We denote with \hypertarget{SK}{$S^k$} a definition concerning conditions on the coefficients (not including the generator) of the BSPDE:
{
\defn{
\label{definition: conditions on coefficients of BSPDE, non-generator}
The condition $S^k$ indicates that $f$ and $g$ are $k$-times continuously differentiable and the partial derivatives up to order $k$ are all in \hyperlink{P}{$\mathcal{P}_T$}.
$G$ and $F$ are $(k + 1)$-times continuously differentiable and the partial derivatives up to order $(k + 1)$ are all uniformly bounded in $(\omega, t, x)$. 
$\varphi$ is $k$-times continuously differentiable, and all partial derivatives of order $0$ to $k$ grow at most polynomially. 
}
}

{
\lemma{
\label{lemma: Rozovskii SPDE result}
Assume \hyperlink{DK}{$D^k$} and \hyperlink{SK}{$S^k$} for some $3 \leq k \in \N{}$. 
Additionally, assume the parabolic condition $2 a - F^* F \succeq 0$ holds. 
Then Eq. \ref{equation: general form of BSPDE} has a unique classical solution $\psi$ in the sense that for every fixed $\omega$ outside a null set, $\psi(\omega, \cdot, \cdot) \in C^{0, k-1}([0, T] \times \R{d}; \R{})$, $\psi$ and its partial derivatives are in \hyperlink{P}{$\mathcal{P}_T(\R{m}; \R{})$}, and $\psi$ solves the integral equation. 
If $\wt{\psi}$ is any other solution of the integral equation, then $\psi$ and $\wt{\psi}$ are indistinguishable. 
If further $f, g$ and $\varphi$ as well as their derivatives up to order $k$ are uniformly bounded in $(\omega, t, x)$, then for any $p > 0$ there exist $q > 0$ and $C > 0$ (only depending on $p$, the dimensions involved, the bounds on $a, b, G$ and $F$, and on $T$), such that for all $| \beta | \leq k - 1$ and $x \in \R{m}$, 
\begin{align*}
\E{ \sup_{t \leq T} | D^\beta \psi(t, x) |^p }
\leq C (1 + |x|^q ) \E{ \norm{\varphi}{k, \infty}^p + \sup_{t \leq T} \norm{ f(t, \cdot) }{k, \infty}^p + \sup_{t \leq T} \norm{ g(t, \cdot) }{k, \infty}^p } .
\end{align*}

\begin{proof}
The lemma is a slight generalization of \cite[Proposition 4.1, p.2302]{Imkeller:2013}, and follows the same argument as the one given there.
\end{proof}
}
}

\subsection{Backward Doubly Stochastic Differential Equations}
\label{subsection: BDSDE}

The theory of backward doubly stochastic differential equations has its origin in the paper by \cite{Pardoux:1994vv}. 
Although it is possible to get a different representation of the solutions of Eq. \ref{equation: general form of BSPDE} by the Method of Stochastic Characteristics \cite{Rozovskii:1990ses}, one benefit of the BDSDE representation is that for fixed $(x, z) \in \R{m} \times \R{n}$, we will have a finite dimensional representation of $\psi(x, z)$ and therefore will be able to apply Gr\"onwall's lemma in the final step of Lemma \ref{lemma: quantitative remainder term}, as part of the main analysis. 

A BDSDE is an integral equation of the form,
\begin{align}
\label{equation: general form of BDSDE}
Y_t = \xi + \int_t^T f(s, \cdot, Y_s, Z_s) ds + \int_t^T \langle g(s, \cdot, Y_s, Z_s) , d \B[s] \rangle - \int_t^T \langle Z_s , dW_s \rangle,
\end{align}
where $f: [0, T] \times \Omega \times \R{} \times \R{n} \rightarrow \R{}$, $g: [0, T] \times \Omega \times \R{} \times \R{n} \rightarrow \R{d}$, and for fixed $y \in \R{}$, $z \in \R{n}$, the processes $(\omega, t) \mapsto f(t, \omega, y, z)$ and $(\omega, t) \mapsto g(t, \omega, y, z)$ are $(\mF^B_{0, T} \vee \mF^W_T) \otimes \mathcal{B}(\R{})$-measurable, and for every $t$, $f(t, \cdot, y, z)$ and $g(t, \cdot, y, z)$ are $\mF_t$-measurable. 
Our definition of $\mF_t$ is, 
\begin{align*}
\mF_t = \mF^B_{t, T} \vee \mF^W_t,
\end{align*}
where $\mF^W_t = \mF^W_{0, t}$. 
Because of this definition, $\mF_t$ is not a filtration; it is neither strictly increasing nor decreasing in $t$. 
Let us now introduce some additional notation for integrability and measurability conditions of the solution of the BDSDEs. 
\hypertarget{H2}{}
{
\defn{
\label{definition: measurable BDSDE condition}
Let $H^2_T(\R{m})$ be the space of measurable $\R{m}$-valued processes $Y$, such that $Y_t$ is $\mF_t$-measurable for almost any $t \in [0, T]$ and 
\begin{align*}
\E{ \int_0^T | Y_t |^2 dt } < \infty. 
\end{align*}
}
}
\hypertarget{S2}{}
{
\defn{
\label{definition: adapted BDSDE condition}
Let $S^2_T(\R{m})$ be the space of continuous adapted $\R{m}$-valued processes $Y$, such that $Y_t$ is $\mF_t$-measurable for every $t \in [0, T]$ and 
\begin{align*}
\E{ \sup_{0 \leq t \leq T} | Y_t |^2 dt } < \infty. 
\end{align*}
}
}

The pair $(Y, Z)$ will be called a solution of Eq. \ref{equation: general form of BDSDE} if $(Y, Z) \in S^2_T(\R{}) \times H^2_T(\R{n})$, and if the pair solves the integral equation. 
We will also write BDSDEs in differential form at times, for example Eq. \ref{equation: general form of BDSDE} in differential form would be, 
\begin{align*}
- d Y_t = f(t, \cdot, Y_t, Z_t) dt + \langle g(t, \cdot, Y_t, Z_t) , d \B \rangle - \langle Z_t , dW_t \rangle.
\end{align*}
With suitable adaptations, all of the following results also hold in the multidimensional case (i.e., $Y \in \R{m}$). 
We restrict to the one dimensional case for simplicity and because ultimately we are only interested in that case.

In \cite{Pardoux:1994vv}, it is shown that under the following conditions, Eq. \ref{equation: general form of BDSDE} has a unique solution: 
\begin{itemize}
\item $\xi \in L^2(\Omega, \mF_T, \Pe; \R{})$,
\item for any $(y, z) \in \R{} \times \R{n}$, we have $f(\cdot , \cdot, y, z) \in$ \hyperlink{H2}{$H^2_T(\R{})$} and $g(\cdot , \cdot, y, z) \in$ \hyperlink{H2}{$H^2_T(\R{d})$}
\item $f$ and $g$ satisfy Lipschitz conditions and $g$ is a contraction in $z$: there exists constants $L > 0$ and $0 < \beta < 1$ such that for any $(\omega, t)$ and $y_1, y_2, z_1, z_2$, 
\begin{align*}
|f(t, \omega, y_1, z_1) - f(t, \omega, y_2, z_2)|^2 
&\leq L (|y_1 - y_2|^2 + |z_1 - z_2|^2), \\
|g(t, \omega, y_1, z_1) - g(t, \omega, y_2, z_2)|^2 
&\leq L|y_1 - y_2|^2 + \beta |z_1 - z_2|^2.
\end{align*}
\end{itemize}
Now we associate a diffusion $X$ to the differential operator $\mathcal{L}$ given in Eq. \ref{equation: general form of BSPDE}. 
To do so, assume \hyperlink{DK}{$D_k$} is satisfied for some $k \geq 2$. 
Then $\sigma \equiv a^{1/2}$ is Lipschitz continuous by \cite[Lemma 2.3.3]{Stroock:2008c}. 
Hence for every $(t, x) \in [0, T] \times \R{m}$, there exists a strong solution of the SDE
\begin{align*}
X^{t, x}_s 
&= x + \int_t^s b(X^{t, x}_s) ds + \int_t^s \sigma(X^{t, x}_s) dW_s, \qquad \text{for} \ s \geq t, \\
X^{t, x}_s
&= x \qquad \text{for} \ s \leq t.
\end{align*}
For the theory of BDSDEs, we must assume that $F$ has the form $F = \alpha \sigma^*$, and here we consider $\alpha \in \R{d \times m}$ a constant matrix. 
We then associate the following BDSDE to Eq. \ref{equation: general form of BSPDE}, 
\begin{align*}
- d Y^{t, x}_s 
&= f(s, X^{t, x}_s) ds + \langle g(s, X^{t, x}_s) + G(s, X^{t, x}_s) Y^{t, x}_s + \alpha Z^{t, x}_s , d \B \rangle - \langle Z^{t, x}_s , dW_s \rangle. \\
Y^{t, x}_T
&= \varphi(X^{t, x}_T).
\end{align*}
Under the assumptions \hyperlink{SK}{$S_k$} and \hyperlink{DK}{$D_k$} for $k \geq 2$, this equation has a unique solution. 
The tuple $(X^{t, x}, Y^{t, x}, Z^{t, x})$ constitutes a forward backward doubly stochastic differential equation (FBDSDE). 

{
\lemma{
\label{lemma: BDSDE implies BSPDE solution}
Assume \hyperlink{SK}{$S_k$} and \hyperlink{DK}{$D_k$} for $k \geq 3$ and $2a - F^* F \succeq 0$. 
Then the unique classical solution $\psi$ of the BSPDE in Eq. \ref{equation: general form of BSPDE} is given by $\psi(t, x) = Y^{t, x}_t$, where $(Y^{t, x}, Z^{t, x})$ is the unique solution of the BDSDE in Eq. \ref{equation: general form of BDSDE}. 
\begin{proof}
See \cite[Theorem 3.1, p.225]{Pardoux:1994vv}.
\end{proof}
}
}

A final remark before we turn to the preliminary estimates and the main analysis where we will use BDSDEs,
we will not be able to get an existence result for classical solutions of the SPDEs in Section \ref{section: main analysis of quantitative convergence} from the theory of BDSDEs. 
This is due to the fact that for this we would need smoothness properties of a square root of $a$. 
But even when $a$ is smooth, in the degenerate elliptic case it does not need to have a smooth square root (see for example \cite[Lemma 2.3.3]{Stroock:2008c}). 
We will instead use the existence result of \cite{Rozovskii:1990ses} stated in Lemma \ref{lemma: Rozovskii SPDE result}, and use the uniqueness result of \cite{Pardoux:1994vv} in our setting. 
This works under Lipschitz continuity of $a^{1/2}$, which we used previously.

\section{Preliminary Estimates}
\label{section: preliminary estimates (chapter: quantitative)}

In this section, we prove several preliminary estimates to be used in Section \ref{section: main analysis of quantitative convergence}. 
We start with results for the moments of the SDE solutions.

\subsection{Estimates on SDE Solutions}

{
\lemma{
\label{lemma: moments on slow process (chapter: quantitative)}
Assume that the drift coefficient $b$, and dispersion coefficient $\sigma$, of the slow motion $\Xe$ are bounded. 
Then for any $p \geq 1$, and every $T > 0$, there exists $C_p > 0$ such that 
\begin{align*}
\sup_{(t, \epsilon) \in [0, T] \times (0, 1] } \E{ | \Xe_t |^p}[(\Xe_0, \Ze_0) = (x, z)] \leq C_p (1 + |x|^p).
\end{align*}
\begin{proof}
The result is trivial since we assume the coefficients to be bounded and consider finite $T$. 
\end{proof}
}
}

{
\lemma{
\label{lemma: moments on fast process (chapter: quantitative)}
Assume $f$ is bounded and that $f$ and $gg^*$ are H\"older continuous in $z$ uniformly in $x$ for some uniform constant.
Assume that the conditions \ref{assumption: positive recurrence} and \ref{assumption: uniform ellipticity} hold. 
Then for any $p > 0$ there exists $C_p > 0$ such that 
\begin{align*}
\sup_{(t, \epsilon, x) \in [0, \infty) \times (0, 1] \times \R{m}} \E{ | \Ze_t |^p}[(\Xe_0, \Ze_0) = (x, z)] \leq C_p (1 + |z|^p). 
\end{align*}

\begin{proof}
The lemma is a slight generalization of a part of \cite[Proposition 5.3, p.2307]{Imkeller:2013}, and the proof follows the same argument as given there.
\end{proof}
}
}

{
\lemma{
\label{lemma: moments of Girsanov (chapter: quantitative)} 
Assume $h$ is bounded, then for $p \geq 1$ and $t \in [0, T]$, 
\begin{align*}
\sup_{\epsilon \in (0, 1]} \sup_{t \leq T} \E*[\Pe] \left| \wt{D}^\epsilon_t \right|^p < \infty .
\end{align*}
\begin{proof}
See the proof of \cite[Lemma 6.5]{Imkeller:2013}.
\end{proof}
}
}

\subsection{Estimates with the Fast Semigroup}

In this section we provide estimates relating to the semigroup of the fast process.

{
\lemma{
\label{lemma: bounds on functions under semigroup (chapter: quantitative)}
Assume \hyperlink{HF}{$HF^{k, l}$}, with $k \in \N{}_0, l \in \N{}$, and let $\theta \in C^{k, j}(\R{m} \times \R{n}; \R{})$ for $j \leq l$ satisfy for some $C, p > 0$
\[
\sum_{|\alpha| \leq k} \sum_{|\beta| \leq j} | D^\alpha_x D^\beta_z \theta(x, z) | \leq C (1 + |x|^p + |z|^p).
\]
Then 
\[
(t, x, z) \mapsto T^{F, x}_t \left( \theta (x, \cdot) \right) (z) \in C^{0, k, j}(\R{+} \times \R{m} \times \R{n}; \R{})
\]
and there exist $C_1, p_1 > 0$, such that for all $(t, x, z) \in [0, \infty) \times \R{m} \times \R{n}$
\[
\sum_{|\alpha| \leq k} \sum_{|\beta| \leq j} | D^\alpha_x D^\beta_z T^{F, x}_t \left( \theta (x, \cdot) \right) (z) | \leq C_1 e^{C_1 t} (1 + |x|^{p_1} + |z|^{p_1}).
\]
If the bound on the derivatives of $\theta$ can be chosen uniformly in $x$, that is, 
\[
\sum_{|\alpha| \leq k} \sum_{|\beta| \leq j} \sup_x | D^\alpha_x D^\beta_z \theta(x, z) | \leq C (1 + |z|^p),
\]
then the bound on the derivatives of $T^{F, x}_t \left( \theta (x, \cdot) \right) (z)$ is also uniform in $x$, 
\[
\sum_{|\alpha| \leq k} \sum_{|\beta| \leq j} \sup_x | D^\alpha_x D^\beta_z T^{F, x}_t \left( \theta (x, \cdot) \right) (z) | \leq C_1 e^{C_1 t} (1 + |z|^{p_1}).
\]

\begin{proof}
The lemma is a slight generalization of \cite[Proposition 5.1]{Imkeller:2013}. 
The proof is the same as in \cite[Proposition 5.1]{Imkeller:2013}.
\end{proof}
}
}

{
\lemma{
\label{lemma: regularity and bound on averaged function (chapter: quantitative)}
Assume \ref{assumption: positive recurrence}, \ref{assumption: uniform ellipticity} and \hyperlink{HF}{$HF^{k, 3}$} for $k \in \N{}_0$. 
Let $\theta \in C^{k, 0}(\R{m} \times \R{n}; \R{})$ satisfy for some $C, p > 0$,
\[
\sum_{|\gamma| \leq k} \sup_x | D^\gamma_x \theta(x, z) | \leq C ( 1 + |z|^p ). 
\]
Then 
\[
x \mapsto \mu_\infty(\theta; x)(x') = \int_{\R{n}} \theta(x', z) \mu_\infty(dz; x) \in C^{k}_b(\R{m}; \R{}).
\]

\begin{proof}
See \cite[Proposition 5.2]{Imkeller:2013}, which contains the same statement. 
\end{proof}
}
}

{
\lemma{
\label{lemma: bounds on centered functions under semigroup (chapter: quantitative)}
Assume \ref{assumption: positive recurrence}, \ref{assumption: uniform ellipticity} and \hyperlink{HF}{$HF^{k, 3}$} with $k \in \N{}_0$. 
Let $\theta \in C^{k, 1}(\R{m} \times \R{n}; \R{})$ satisfy the growth condition, 
\[
\sum_{|\alpha| \leq k} \sum_{|\beta| \leq 1} \sup_x | D^\alpha_x D^\beta_z \theta(x, z) | 
\leq C (1 + |z|^p),
\]
for some $C, p > 0$ .
Assume additionally that $\theta$ satisfies the centering condition, 
\[
\int_{\R{n}} \theta(x, z) \mu_\infty(dz; x) = 0, \quad \forall x \in \R{m}.
\]
Then 
\[
(x, z) \mapsto \int_0^\infty T^{F, x}_t(\theta(x, \cdot))(z) dt \in C^{k, 1}(\R{m} \times \R{n}; \R{}),
\]
and for every $q > 0$ there exists $C', q' > 0$, such that, 
\[
\sum_{|\alpha| \leq k} \sum_{|\beta| \leq 1} \int_0^\infty \sup_x | D^\alpha_x D^\beta_z T^{F, x}_t(\theta(x, \cdot))(z) |^q dt 
\leq C' (1 + |z|^{q'}).
\]

\begin{proof}
The lemma is a slight generalization of a part of \cite[Proposition 5.2]{Imkeller:2013}, and the proof follows the same argument as given there.
\end{proof}
}
}

{
\lemma{
\label{lemma: regularity and integrability of simple averaged coefficients (chapter: quantitative)}
Assume \ref{assumption: positive recurrence}, \ref{assumption: uniform ellipticity} and \hyperlink{HF}{$HF^{k, 3}$} with $k \in \N{}_0$.
If \hyperlink{HS}{$HS^{k, 1}$} holds, then $\overbar{b}, \overbar{\sigma}, \overbar{a} \in C^k_b$.
Similarly, if \hyperlink{HO}{$HO^{k, 1}$} holds, then $\overbar{h} \in C^k_b$. 

\begin{proof}
The result follows from Lemma \ref{lemma: bounds on centered functions under semigroup (chapter: quantitative)}.
\end{proof}
}
}

\subsection{Estimates for the Corrector Term}

We now introduce a few lemmas that will help to streamline the main ideas in the analysis of the corrector term given in Lemma \ref{lemma: quantitative corrector term}.
{
\lemma{
\label{lemma: corrector analysis, time integral}
Assume that $u \in C([0, T] \times \R{m}; \R{})$ and is an element of \hyperlink{P}{$\mathcal{P}_T(\R{m}; \R{})$}, 
that the conditions \ref{assumption: positive recurrence}, \ref{assumption: uniform ellipticity}, and \hyperlink{HF}{$HF^{0, 3}$} hold, and that $\psi \in C^{0, 1}_b(\R{m} \times \R{n}; \R{})$. 
Assume additionally that $\psi$ satisfies the centering condition, 
\[
\int_{\R{n}} \psi(x, z) \mu_\infty(dz; x) = 0, \quad \forall x \in \R{m}.
\]
Then given $(x, z) \in \R{m} \times \R{n}$ and $t \in [0, T]$, there exists $q > 0$ such that
\begin{align*}
\E{ \left| \E{ \int_t^T \psi(x, \Zextz) u_s(x) ds }[\F] \right|^p }
\lesssim \epsilon^{2 p} ( 1 + |z|^q ) \E{ \sup_{t \leq s \leq T} \left| u_s(x) \right|^p }.
\end{align*}
Here $\Zex$ is the diffusion process with generator $\frac{1}{\epsilon^2} \mG_F$ \emph{(}in  particular, the Brownian motion driving $\Zex$ is independent of the Brownian motion $\B$ that generates the backward filtration $\F$\emph{)}.

\begin{proof}
Because $u_s$ is measurable with respect to $\F[s, T]$ and $\Zextz[]$ is independent of $B$, we get from the conditional expectation with respect to $\F$ and definition of the semigroup $T^{F, x}$ the following identity
\begin{align*}
\E{ \int_t^T \psi(x, \Zextz) u_s(x) ds }[\F] 
= \int_t^T \E{ \psi(x, \Zextz) } u_s(x) ds 
= \int_t^T T^{F, x}_{(s - t) / \epsilon^2}(\psi(x, \cdot))(z) u_s(x) ds.
\end{align*}
Now taking the absolute value, using H\"older's inequality, and removing $| u_s(x) |$ from the integral by taking the supremum over $[t, T]$ gives
\begin{align}
\label{equation: corrector analysis, time integral, time change step}
\left| \int_t^T T^{F, x}_{(s - t) / \epsilon^2}(\psi(x, \cdot))(z) u_s(x) ds \right| 
\lesssim \sup_{t \leq s \leq T} | u_s(x) | \int_t^T \left| T^{F, x}_{(s - t) / \epsilon^2}(\psi(x, \cdot))(z) \right| ds.
\end{align}
Performing a time reparametrization and then using Lemma \ref{lemma: bounds on centered functions under semigroup (chapter: quantitative)}, we have for some $q' > 0$, 
\begin{align*}
\sup_{t \leq s \leq T} | u_s(x) | \int_t^T \left| T^{F, x}_{(s - t) / \epsilon^2}(\psi(x, \cdot))(z) \right| ds
&\lesssim \epsilon^2 \sup_{t \leq s \leq T} | u_s(x) | \int_0^{(T - t) / \epsilon^2} \left| T^{F, x}_r(\psi(x, \cdot))(z) \right| dr \\
&\lesssim \epsilon^2 \sup_{t \leq s \leq T} | u_s(x) | \int_0^\infty \left| T^{F, x}_r(\psi(x, \cdot))(z) \right| dr \\
&\lesssim \epsilon^2 (1 + |z|^{q'}) \sup_{t \leq s \leq T} | u_s(x) | .
\end{align*}
Lastly, taking the $p$-th power and applying the expectation gives the desired result.
\end{proof}
}
}

{
\lemma{
\label{lemma: corrector analysis, stochastic integral}
Assume that $u \in C([0, T] \times \R{m}; \R{k})$ and is an element of \hyperlink{P}{$\mathcal{P}_T(\R{m}; \R{k})$} for $k \geq 1$.
Assume \ref{assumption: positive recurrence}, \ref{assumption: uniform ellipticity} and \hyperlink{HF}{$HF^{0, 3}$}. 
Let $\psi \in C^{0, 1}_b(\R{m} \times \R{n}; \R{d \times k})$. 
Assume additionally that $\psi$ satisfies the centering condition, 
\[
\int_{\R{n}} \psi(x, z) \mu_\infty(dz; x) = 0, \quad \forall x \in \R{m}.
\] 
Then given $(x, z) \in \R{m} \times \R{n}$ and $t \in [0, T]$, there exists $q > 0$ such that
\begin{align*}
\E{ \left| \E{ \int_t^T  \langle \psi(x, \Zextz) u_s(x) , d \B[s] \rangle }[\F] \right|^p }
\lesssim \epsilon^p ( 1 + |z|^q ) \E{ \sup_{t \leq s \leq T} \left| u_s(x) \right|^p }.
\end{align*}
Here $\Zex$ is the diffusion process with generator $\frac{1}{\epsilon^2} \mG_F$ \emph{(}in  particular, the Brownian motion driving $\Zex$ is independent of the Brownian motion $\B$ that generates the backward filtration $\F$\emph{)}.

\begin{proof}
Because $u_s$ is measurable with respect to $\F[s, T]$ and $\Zextz[]$ is independent of $B$, we get from the conditional expectation with respect to $\F$ and definition of the semigroup $T^{F, x}$ the following identity
\begin{align*}
\E{ \int_t^T  \langle \psi(x, \Zextz) u_s(x) , d \B[s] \rangle }[\F]
&= \int_t^T \langle \E{ \psi(x, \Zextz) } u_s(x) , d \B[s] \rangle \\
&= \int_t^T \langle T^{F, x}_{(s - t) / \epsilon^2}(\psi(x, \cdot))(z) u_s(x) , d \B[s] \rangle .
\end{align*}
Now by application of the Burkholder-Davis-Gundy inequality we get
\begin{align*}
\E{ \left| \int_t^T \langle T^{F, x}_{(s - t) / \epsilon^2}(\psi(x, \cdot))(z) u_s(x) , d \B[s] \rangle \right|^p }
&\lesssim \E{ \left\langle \int_t^T \langle T^{F, x}_{(s - t) / \epsilon^2}(\psi(x, \cdot))(z) u_s(x) , d \B[s] \rangle \right\rangle^{p / 2} }.
\end{align*}
Computing the quadratic variation gives 
\begin{align}
\label{equation: quadratic variation in proposition for corrector analysis}
\left\langle \int_t^T \langle T^{F, x}_{(s - t) / \epsilon^2}(\psi(x, \cdot))(z) u_s(x) , d \B[s] \rangle \right\rangle
&= \int_t^T | T^{F, x}_{(s - t) / \epsilon^2}(\psi(x, \cdot))(z) u_s(x) |^2 ds .
\end{align}
In the case that $u$ is real-valued, the integrand for the right side of Eq. \ref{equation: quadratic variation in proposition for corrector analysis} is bounded by 
\begin{align}
\label{equation: inequality in proposition for corrector analysis}
| T^{F, x}_{(s - t) / \epsilon^2}(\psi(x, \cdot))(z) u_s(x) |^2
\leq | u_s(x) |^2 | T^{F, x}_{(s - t) / \epsilon^2}(\psi(x, \cdot))(z) |^2.
\end{align}
Similarly, in the case where $u$ is an  $\R{k}$-valued vector for some $k > 1$, then $T^{F, x}_{(s - t) / \epsilon^2}(\psi(x, \cdot))(z)$ takes values in $\R{d \times k}$, and letting 
$A_s  \equiv T^{F, x}_{(s - t) / \epsilon^2}(\psi(x, \cdot))(z)$ for brevity, we have $| A_s u_s(x) |^2 \lesssim | u_s(x) |^2 \op{Tr}( A^*_s A_s ) = | u_s(x) |^2 | A_s |^2$.
Therefore we have the same inequality for the integrand.

Therefore using Eq. \ref{equation: inequality in proposition for corrector analysis} in the quadratic variation of Eq. \ref{equation: quadratic variation in proposition for corrector analysis} and then taking the function $| u_s(x) |$ outside the integral by using its supremum value over $[t, T]$, we get
\begin{align*}
\left\langle \int_t^T \langle T^{F, x}_{(s - t) / \epsilon^2}(\psi(x, \cdot))(z) u_s(x) , d \B[s] \rangle \right\rangle
&\leq \int_t^T | T^{F, x}_{(s - t) / \epsilon^2}(\psi(x, \cdot))(z) |^2 | u_s(x) |^2 ds \\
&= \sup_{t \leq s \leq T} | u_s(x) |^2 \int_t^T | T^{F, x}_{(s - t) / \epsilon^2}(\psi(x, \cdot))(z) |^2  ds .
\end{align*}
We now perform a time reparametrization and use Lemma \ref{lemma: bounds on centered functions under semigroup (chapter: quantitative)}, so that for some $q' > 0$ we get
\begin{align}
\sup_{t \leq s \leq T} | u_s(x) |^2 \int_t^T | T^{F, x}_{(s - t) / \epsilon^2}(\psi(x, \cdot))(z) |^2  ds \nonumber
&= \epsilon^2 \sup_{t \leq s \leq T} | u_s(x) |^2 \int_0^{(T - t) / \epsilon^2} | T^{F, x}_r (\psi(x, \cdot))(z) |^2  dr \nonumber \\
&= \epsilon^2 \sup_{t \leq s \leq T} | u_s(x) |^2 \int_0^\infty | T^{F, x}_r (\psi(x, \cdot))(z) |^2  dr \nonumber \\
&\lesssim \epsilon^2 (1 + |z|^{q'}) \sup_{t \leq s \leq T} | u_s(x) |^2 . \label{equation: final term in corrector analysis} 
\end{align}
To see that the last step still holds in the case that $\psi$ is matrix-valued, consider the following relations
\begin{align*}
| T^{F, x}_r (\psi(x, \cdot))(z) |^2
&= \sum_{i = 1}^d \left( T^{F, x}_r (\psi(x, \cdot))(z) T^{F, x}_r (\psi(x, \cdot))(z)^* \right)_{ii} \\
&= \sum_{i, j = 1}^d \left( T^{F, x}_r (\psi(x, \cdot))(z) \right)_{ij}^2 
= \sum_{i, j = 1}^d | T^{F, x}_r (\psi(x, \cdot)_{ij} )(z) |, 
\end{align*}
which shows that the same analysis holds, but now for a summation over the centered entries of $\psi$. 
Finally, taking the $p/2$ power and applying the expectation to Eq. \ref{equation: final term in corrector analysis} gives the desired result.
\end{proof}
}
}

\section{Main Analysis}
\label{section: main analysis of quantitative convergence}

\subsection{Moment Estimates for Dual Processes}

In this section, we compute the main estimates for $v^0, \psi$ and $R$ associated with an arbitrary fixed test function $\varphi \in C^2_b$. 
The estimates for $\psi$ and $R$ are then used in Section \ref{subsection: estimates of dual and filter error} to prove Theorem \ref{theorem: main result of quantitative convergence}.

{
\lem{
\label{lemma: quantitative averaged term}
Let $3 \leq k \in \N{}$ and assume $\overbar{b}, \overbar{a}, \varphi \in C^k_b$ and $\overbar{h}, \overbar{\sigma} \in C^{k + 1}_b$. 
Then $v^0 \in C^{0, k - 1}([0, T] \times \R{m}; \R{})$, and for any $p \geq 1$ there exist $q > 0$, such that for all $x \in \R{m}$, 
\begin{align*}
\sum_{|j| \leq k - 1} \E{\sup_{0 \leq t \leq T} | D^j_x v^0_t (x) |^p} 
\lesssim ( 1 + |x|^q ) \norm{\varphi}{k, \infty}^p.
\end{align*}
In particular, $v^0$ and all its partial derivatives up to order $(0, k - 1)$ are in \hyperlink{P}{$\mathcal{P}_T(\R{m}; \R{})$}. 
\begin{proof} 
The result follows from Lemma \ref{lemma: Rozovskii SPDE result}. 
The only condition from Lemma \ref{lemma: Rozovskii SPDE result} that is not immediately obvious is the parabolic condition, $2\overbar{a} - \overbar{\sigma} \alpha^* \alpha \overbar{\sigma}^* \succeq 0$. 
This condition indeed holds for the same reason as given in Section \ref{section: the averaged conditional distributions} and the fact that $I - \alpha^* \alpha \succeq 0$, where $I$ is the identity matrix (recall that $\alpha$ was redefined in Section \ref{section: problem statement} as $\alpha \leftarrow \kappa^{-1} \alpha$, and note that $(\kappa^{-1})^* = (\kappa^*)^{-1}$). 
\end{proof}
}
}

{
\lem{
\label{lemma: quantitative corrector term}
Let $3 \leq k, l \in \N{}$ and assume \ref{assumption: positive recurrence}, \ref{assumption: uniform ellipticity}, \hyperlink{HF}{$HF^{k, l}$}, \hyperlink{HS}{$HS^{k, l}$}, \hyperlink{HO}{$HO^{k, l}$}, and that $\overbar{\sigma}, \overbar{a}, \overbar{b}, \overbar{h} \in C^k_b$. 
Let $v^0 \in C^{0, k}([0, T] \times \R{m}; \R{})$, and assume that all its partial derivatives in $x$ up to order $k$ are in \hyperlink{P}{$\mathcal{P}_T(\R{m}; \R{})$}. 

Then $\psi \in C^{0, k - 1, l - 1}([0, T] \times \R{m} \times \R{n}; \R{})$, and $\psi$ as well as its partial derivatives up to order $(0, k - 1, l - 1)$ are in \hyperlink{P}{$\mathcal{P}_T(\R{m} \times \R{n}; \R{})$}. 
For any $p \geq 1$ there exists $q > 0$, such that for any $(x, z) \in \R{m} \times  \R{n}$ and any $\epsilon \in (0, 1)$
\begin{align*}
\sum_{|\beta| \leq k - 2} \sup_{0 \leq t \leq T} \E{|D^\beta_x \psi_t(x, z)|^p} 
\lesssim \epsilon^p ( 1 + |z|^q ) \sum_{ |\beta| \leq k} \E{\sup_{0 \leq t \leq T} | D^\beta_x v^0_t(x)|^p}.
\end{align*}

\begin{proof}
$\psi_t(x, z)$ solves the following BSPDE
\EquationAligned*{
-d \psi_t (x, z)
&= \left[ \frac{1}{\epsilon^2} \mG_F \psi_t (x, z) + (\mG_S - \overbar{\mG_S}) v^0_t (x, z) \right] dt + \langle v^0_t (h - \overbar{h}) (x, z) , d \B \rangle + \langle \alpha (\sigma - \overbar{\sigma})^* \nabla_x v^0_t (x, z) , d \B \rangle, \\
\psi_T(x, z) 
&= 0 .
}
Existence of the solution $\psi$ and its derivatives as well as the polynomial growth follow from Lemma \ref{lemma: Rozovskii SPDE result}. 
From Lemma \ref{lemma: BDSDE implies BSPDE solution}, the solution, $\psi_t(x, z)$, has a representation in terms of a FBDSDE, $\psi_t(x, z) = \theta_t^{t, x, z}$.  
Where $\theta$ is a component of the pair of processes $(\theta, \gamma^\epsilon)$ satisfying the BDSDE
\EquationAligned{
\label{equation: corrector term, BDSDE}
-d \theta_s^{t, x, z} 
&= \left[ \mG_S( x, \Zextz)  - \overbar{\mG_S}(x) \right] v^0_s(x) ds \\
&+ \langle v^0_s(x) (h( x, \Zextz) - \overbar{h}(x)) , d \B[s] \rangle
+ \langle \alpha (\sigma (x, \Zextz) - \overbar{\sigma}(x))^* \nabla_x v^0_s(x) , d \B[s] \rangle \\
&- \langle \gamma_s^{\epsilon; (t, x, z)} , dV_s \rangle , \\
\theta_T^{t, x, z} &= 0,
}
and $(x, \Zextz)$ is a joint diffusion process with $\Xetx$ having the zero generator, 
\[
\Xetx = x, \qquad \forall \, s \in [t, T],
\]
and $\Zextz$ satisfying the stochastic differential equation
\EquationAligned*{
d\Zextz 
&= \frac{1}{\epsilon^2} f(x, \Zextz) ds + \frac{1}{\epsilon} g(x, \Zextz) dV_s, \qquad s \geq t, \\
\Zextz &= \ds z, \qquad s \leq t.
}
The second component of the pair $(\theta^{t, x, z}_t, \gamma^{\epsilon; (t, x, z)}_t)$, has a representation as 
\[
\gamma^{\epsilon; (t, x, z)}_t = \frac{1}{ \epsilon} g^* \nabla_z \psi_t(x, z).
\]
For brevity, let us temporarily drop from the notation,  superscripts and part of superscripts that indicate initial conditions (for example, $(t, x, z)$ and $(t, z)$). 

Since $\psi_t$ is $\F$-measurable, so is $\theta_t$, and therefore conditioning $\theta_t$ on $\F$ gives $\theta_t = \E{\theta_t}[\F]$.  
We also observe that $V$ and $B$ are independent. 
Therefore, $V$ is a Brownian motion in the larger filtration $(\mF^V_s \vee \F)_{s \in [0, T]}$.
Hence, from an application of the tower property of conditional expectation, we have
\begin{align*}
\E{\int_t^T \langle \gamma^\epsilon_s , dV_s \rangle }[\F]
= \E{ \E{ \int_t^T \langle \gamma^\epsilon_s , dV_s \rangle }[\mF^V_t \vee \F] }[\F]
= 0,
\end{align*}
and therefore, 
\EquationAligned{
\label{equation: corrector term, BDSDE, application of conditional expectation}
\theta_t 
&= \E{ \int_t^T \left[ \mG_S( x, \Zex_s)  - \overbar{\mG_S}(x) \right] v^0_s(x) ds}[\F]  \\
&+ \E{ \int_t^T \langle v^0_s(x) (h( x, \Zex_s) - \overbar{h}(x)) , d \B[s] \rangle }[\F] \\
&+ \E{ \int_t^T \langle \alpha (\sigma( x, \Zex_s) - \overbar{\sigma}(x))^* \nabla_x v^0_s(x) , d \B[s] \rangle }[\F]. 
}
The $p$-th moment is therefore bounded as follows, 
\begin{align}
\E{ \left| \theta_t \right|^p }
&\lesssim \E{ \left| \E{ \int_t^T \left[ \mG_S( x, \Zex_s)  - \overbar{\mG_S}(x) \right] v^0_s(x) ds}[\F] \right|^p } \label{equation: corrector term, BDSDE, difference of generator} \\
&+ \E{ \left| \E{ \int_t^T \langle v^0_s(x) (h( x, \Zex_s) - \overbar{h}(x)) , d \B[s] \rangle }[\F] \right|^p } \label{equation: corrector term, BDSDE, difference of sensor} \\
&+ \E{ \left| \E{ \int_t^T \langle \alpha (\sigma( x, \Zex_s) - \overbar{\sigma}(x))^* \nabla_x v^0_s(x) , d \B[s] \rangle }[\F] \right|^p }. \label{equation: corrector term, BDSDE, difference of correlated}
\end{align}

The first term on the right side of Eq. \ref{equation: corrector term, BDSDE, difference of generator} has an integrand that can be written as
\begin{align*}
\left[ \mG_S( x, \Zex_s)  - \overbar{\mG_S}(x) \right] v^0_s(x)
= \sum_{i = 1}^m (b - \overbar{b})_i \frac{\partial}{\partial x_i} v^0_s (x, \Zex_s) + \frac{1}{2} \sum_{i, j = 1}^m (a - \overbar{a})_{ij} \frac{\partial^2}{\partial x_i \partial x_j} v^0_s (x, \Zex_s),
\end{align*}
which shows that this term is a summation of terms that fit the conditions of Lemma \ref{lemma: corrector analysis, time integral} (i.e., a centered function driven by $\Zex$ and multiplied with a term that has the correct bounds and measurability properties) and therefore we get for some $q_0 > 0$ the following estimate for this term 
\begin{align}
\label{equation: corrector term, BDSDE, estimate of difference of generator}
\E{ \left| \E{ \int_t^T \left[ \mG_S( x, \Zex_s)  - \overbar{\mG_S}(x) \right] v^0_s(x) ds}[\F] \right|^p } 
\lesssim \epsilon^{2 p} ( 1 + |z|^{q_0} ) \sum_{1 \leq |\beta| \leq 2} \E{ \sup_{t \leq s \leq T} \left| D^\beta_x v^0_s(x) \right|^p }.
\end{align}

The second term, Eq. \ref{equation: corrector term, BDSDE, difference of sensor}, fits the assumptions of Lemma \ref{lemma: corrector analysis, stochastic integral}, and therefore we get for some $q_1 > 0$ the following estimate for this term
\begin{align}
\label{equation: corrector term, BDSDE, estimate of difference of sensor}
\E{ \left| \E{ \int_t^T \langle v^0_s(x) (h( x, \Zex_s) - \overbar{h}(x)) , d \B[s] \rangle }[\F] \right|^p }
\lesssim \epsilon^p ( 1 + |z|^{q_1} ) \E{ \sup_{t \leq s \leq T} \left| v^0_s(x) \right|^p }.
\end{align}
Unlike the time integral term, we only get $\epsilon^p$ for this estimate because of the application of the Burkholder-Davis-Gundy inequality in the proof of Lemma \ref{lemma: corrector analysis, stochastic integral}.

Lemma \ref{lemma: corrector analysis, stochastic integral} also covers the case where the integrand of the stochastic integral is a matrix-vector product, as occurs in Eq. \ref{equation: corrector term, BDSDE, difference of correlated}.
Because each entry in Eq. \ref{equation: corrector term, BDSDE, difference of correlated} is centered and $v^0$ meets the required conditions of Lemma \ref{lemma: corrector analysis, stochastic integral}, we get for some $q_2 > 0$ the following estimate 
\begin{multline}
\label{equation: corrector term, BDSDE, estimate of difference of correlated}
\E{ \left| \E{ \int_t^T \langle \alpha (\sigma( x, \Zex_s) - \overbar{\sigma}(x))^* \nabla_x v^0_s(x) , d \B[s] \rangle }[\F] \right|^p } \\
\lesssim \epsilon^p ( 1 + |z|^{q_2} ) \sum_{|\beta| = 1} \E{ \sup_{t \leq s \leq T} \left| D^\beta_x v^0_s(x) \right|^p } .
\end{multline}

Collecting the estimates from Eqs. \ref{equation: corrector term, BDSDE, estimate of difference of generator}, \ref{equation: corrector term, BDSDE, estimate of difference of sensor}, and \ref{equation: corrector term, BDSDE, estimate of difference of correlated}, we get for some $q_3 > 0$ the following estimate for the BDSDE solution, 
\begin{align}
\label{equation: corrector term, BDSDE, estimate}
\E{ \left| \theta_t \right|^p }
\lesssim \epsilon^p ( 1 + |z|^{q_3} ) \sum_{|\beta| \leq 2} \E{ \sup_{t \leq s \leq T} \left| D^\beta_x v^0_s(x) \right|^p }.
\end{align}

We will also need estimates of the first and second-order derivatives of $\psi$ in the $x$-component for estimating the remainder term $R$ in Lemma \ref{lemma: quantitative remainder term} (see for instance Eq. \ref{equation: BSPDE, remainder term}). 
Therefore consider taking a first-order partial derivative of Eq. \ref{equation: corrector term, BDSDE, application of conditional expectation}, and then taking the $p$-th moment, and separating terms on the right side of the equation by H\"older's inequality, 
\begin{align}
\E{ \left| \frac{\partial}{\partial x_k} \theta_t \right|^p }
&\lesssim \E{ \left| \frac{\partial}{\partial x_k} \E{ \int_t^T \left[ \mG_S( x, \Zex_s)  - \overbar{\mG_S}(x) \right] v^0_s(x) ds}[\F] \right|^p } \label{equation: corrector term, BDSDE-derivative, difference of generator} \\
&+ \E{ \left| \frac{\partial}{\partial x_k} \E{ \int_t^T \langle v^0_s(x) (h( x, \Zex_s) - \overbar{h}(x)) , d \B[s] \rangle }[\F] \right|^p } \label{equation: corrector term, BDSDE-derivative, difference of sensor} \\
&+ \E{ \left| \frac{\partial}{\partial x_k} \E{ \int_t^T \langle \alpha (\sigma( x, \Zex_s) - \overbar{\sigma}(x))^* \nabla_x v^0_s(x) , d \B[s] \rangle }[\F] \right|^p }. \label{equation: corrector term, BDSDE-derivative, difference of correlated}
\end{align}

Just as we dealt with Eq. \ref{equation: corrector term, BDSDE, difference of generator} by first expanding the difference of the generators we can do the same for the right side of Eq. \ref{equation: corrector term, BDSDE-derivative, difference of generator} to get 
\begin{multline*}
\E{ \left| \frac{\partial}{\partial x_k} \E{ \int_t^T \left[ \mG_S( x, \Zex_s)  - \overbar{\mG_S}(x) \right] v^0_s(x) ds}[\F] \right|^p } \\
\lesssim \sum_{1 \leq |\beta| \leq 2} \E{ \left| \frac{\partial}{\partial x_k} \E{ \int_t^T \psi^\beta ( x, \Zex_s) D^\beta_x  v^0_s(x) ds}[\F] \right|^p },
\end{multline*}
where $\psi^\beta$ is either an entry of $b - \overbar{b}$ or $\frac{1}{2} (a - \overbar{a})$, and hence centered. 
Now following the same arguments as in the proof of Lemma \ref{lemma: corrector analysis, time integral}, we are able to get for any multiindex $1 \leq |\beta| \leq 2$, 
\begin{align*}
\E{ \left| \frac{\partial}{\partial x_k} \E{ \int_t^T \psi^\beta ( x, \Zex_s) D^\beta_x  v^0_s(x) ds}[\F] \right|^p }
&= \E{ \left| \frac{\partial}{\partial x_k} \int_t^T \E{ \psi^\beta ( x, \Zex_s) } D^\beta_x  v^0_s(x) ds \right|^p } \\
&= \E{ \left| \frac{\partial}{\partial x_k} \int_t^T T^{F, x}_{(s - t) / \epsilon^2}( \psi^\beta (x, \cdot))(z) D^\beta_x  v^0_s(x) ds \right|^p }.
\end{align*}
Distributing the derivative inside the time integral now gives (ignoring the $p$-th power and expectation for clarity in the next argument)
\begin{align*}
\left| \frac{\partial}{\partial x_k} \int_t^T T^{F, x}_{(s - t) / \epsilon^2}( \psi^\beta (x, \cdot))(z) D^\beta_x  v^0_s(x) ds \right| 
&\leq \left| \int_t^T \frac{\partial}{\partial x_k} T^{F, x}_{(s - t) / \epsilon^2}( \psi^\beta (x, \cdot))(z) D^\beta_x  v^0_s(x) ds \right| \\
&+ \left| \int_t^T T^{F, x}_{(s - t) / \epsilon^2}( \psi^\beta (x, \cdot))(z) \frac{\partial}{\partial x_k} D^\beta_x  v^0_s(x) ds \right| .
\end{align*}
Estimates for both terms are now achieved by applying the procedure in the proof of Lemma \ref{lemma: corrector analysis, time integral} starting from Eq. \ref{equation: corrector analysis, time integral, time change step} onwards (and using the fact that Lemma \ref{lemma: bounds on centered functions under semigroup (chapter: quantitative)} gives bounds for the derivative of the semigroup) to get for some $q_4 > 0$
\begin{multline}
\label{equation: corrector term, BDSDE-derivative, estimate of difference of generator}
\E{ \left| \frac{\partial}{\partial x_k} \E{ \int_t^T \left[ \mG_S( x, \Zex_s)  - \overbar{\mG_S}(x) \right] v^0_s(x) ds}[\F] \right|^p } \\
\lesssim \epsilon^{2 p} ( 1 + |z|^{q_4} ) \sum_{1 \leq |\beta| \leq 3} \E{ \sup_{t \leq s \leq T} \left| D^\beta_x v^0_s(x) \right|^p }.
\end{multline}

Turning our attention now to Eq. \ref{equation: corrector term, BDSDE-derivative, difference of sensor}, we follow the procedure of Lemma \ref{lemma: corrector analysis, stochastic integral} to interchange the conditional expectation and stochastic integration, and then because of \hyperlink{HO}{$HO^{k + 1, l + 1}$}, we can interchange ordinary differentiation and stochastic integration \cite{Karandikar:1983}, and distribute the derivative to get
\begin{align*}
\E{ \left| \frac{\partial}{\partial x_k} \E{ \int_t^T \langle v^0_s(x) (h( x, \Zex_s) - \overbar{h}(x)) , d \B[s] \rangle }[\F] \right|^p } \\
&\hspace{-30mm}\lesssim \E{ \left| \int_t^T \langle \frac{\partial}{\partial x_k} v^0_s(x) T^{F, x}_{(s - t) / \epsilon^2}( (h - \overbar{h})(x, \cdot))(z) , d \B[s] \rangle \right|^p } \\
&\hspace{-30mm}+ \E{ \left| \int_t^T \langle v^0_s(x) \frac{\partial}{\partial x_k} T^{F, x}_{(s - t) / \epsilon^2}( (h - \overbar{h})(x, \cdot))(z) , d \B[s] \rangle \right|^p }. 
\end{align*}
Estimates for both terms on the right side of the equation now follow from the argument in the proof of Lemma \ref{lemma: corrector analysis, stochastic integral}, starting from the application of the Burkholder-Davis Gundy inequality (and again using the fact that Lemma \ref{lemma: bounds on centered functions under semigroup (chapter: quantitative)} gives bounds for the derivative of the semigroup), to yield for some $q_5 > 0$
\begin{multline}
\label{equation: corrector term, BDSDE-derivative, estimate of difference of sensor}
\E{ \left| \frac{\partial}{\partial x_k} \E{ \int_t^T \langle v^0_s(x) (h( x, \Zex_s) - \overbar{h}(x)) , d \B[s] \rangle }[\F] \right|^p } \\
\lesssim \epsilon^p ( 1 + |z|^{q_5} ) \sum_{|\beta| \leq 1} \E{ \sup_{t \leq s \leq T} \left| D^\beta_x v^0_s(x) \right|^p }.
\end{multline}

The last term to address is Eq. \ref{equation: corrector term, BDSDE-derivative, difference of correlated}.
Just as we did when handling Eq. \ref{equation: corrector term, BDSDE-derivative, difference of sensor}, we follow the procedure of Lemma \ref{lemma: corrector analysis, stochastic integral} to interchange the conditional expectation and stochastic integration, and then because of \hyperlink{HO}{$HO^{k + 1, l + 1}$} interchange ordinary differentiation and stochastic integration \cite{Karandikar:1983}, and distribute the derivative to get
\begin{align*}
\E{ \left| \frac{\partial}{\partial x_k} \E{ \int_t^T \langle \alpha (\sigma( x, \Zex_s) - \overbar{\sigma}(x))^* \nabla_x v^0_s(x) , d \B[s] \rangle }[\F] \right|^p } \\
&\hspace{-30mm}\lesssim \E{ \left| \int_t^T \langle \frac{\partial}{\partial x_k} T^{F, x}_{(s - t) / \epsilon^2}( \alpha (\sigma - \overbar{\sigma})^* (x, \cdot) )(z) \nabla_x v^0_s(x) , d \B[s] \rangle  \right|^p } \\
&\hspace{-30mm}+ \E{ \left| \int_t^T \langle T^{F, x}_{(s - t) / \epsilon^2}( \alpha (\sigma - \overbar{\sigma})^* (x, \cdot) )(z) \frac{\partial}{\partial x_k} \nabla_x v^0_s(x) , d \B[s] \rangle \right|^p }. 
\end{align*}
Estimates for both terms on the right side of the equation now follow from the argument in the proof of Lemma \ref{lemma: corrector analysis, stochastic integral}, starting from the application of the Burkholder-Davis-Gundy inequality (and again using the fact that Lemma \ref{lemma: bounds on centered functions under semigroup (chapter: quantitative)} gives bounds for the derivative of the semigroup), to yield for some $q_6 > 0$
\begin{multline}
\label{equation: corrector term, BDSDE-derivative, estimate of difference of correlated}
\E{ \left| \frac{\partial}{\partial x_k} \E{ \int_t^T \langle \alpha (\sigma( x, \Zex_s) - \overbar{\sigma}(x))^* \nabla_x v^0_s(x) , d \B[s] \rangle }[\F] \right|^p } \\
\lesssim \epsilon^p ( 1 + |z|^{q_6} ) \sum_{1 \leq |\beta| \leq 2} \E{ \sup_{t \leq s \leq T} \left| D^\beta_x v^0_s(x) \right|^p }.
\end{multline}

Collecting the estimates from Eqs. \ref{equation: corrector term, BDSDE-derivative, estimate of difference of generator}, \ref{equation: corrector term, BDSDE-derivative, estimate of difference of sensor}, and \ref{equation: corrector term, BDSDE-derivative, estimate of difference of correlated} then yields for some $q_7 > 0$
\begin{align*}
\E{ \left| \frac{\partial}{\partial x_k} \theta_t \right|^p }
\lesssim \epsilon^p ( 1 + |z|^{q_7} ) \sum_{|\beta| \leq 3} \E{ \sup_{t \leq s \leq T} \left| D^\beta_x v^0_s(x) \right|^p }.
\end{align*}

The procedure to take higher-order derivatives is the same as that for the first-order derivatives (simply involving more terms), and therefore taking the supremum of the estimates of these derivatives over $[0, T]$ and summing the terms, we get for some $q > 0$
\begin{align*}
\sum_{|\beta| \leq k - 1} \sup_{0 \leq t \leq T} \E{ \left| D^\beta_x \theta_t \right|^p }
\lesssim \epsilon^p ( 1 + |z|^q ) \sum_{|\beta| \leq k + 1} \E{ \sup_{0 \leq t \leq T} \left| D^\beta_x v^0_t(x) \right|^p }.
\end{align*}
\end{proof}
}
}

{
\lem{
\label{lemma: quantitative remainder term}
Let $3 \leq k, l \in \N{}$ and assume \ref{assumption: regularity and boundedness of fast coefficients}, \ref{assumption: regularity and boundedness of slow coefficients}, $\sigma \in C^{k + 1, l + 1}_b$, and 
\hyperlink{HO}{$HO^{k + 1, l + 1}$}.
Let $\psi \in C^{0, k + 2, l}([0, T] \times \R{m} \times \R{n}; \R{})$ and assume that all its partial derivatives up to order $(0, k + 2, l)$ are in 
\hyperlink{P}{$\mathcal{P}_T(\R{m} \times \R{n}; \R{})$}. 
Then for any $p > 2$,
we have that for any $(x, z) \in \R{m} \times \R{n}$, $\epsilon \in (0, 1)$, and $t \in [0, T]$, 
\begin{align*}
\E{|R_t(x, z)|^p} 
\lesssim \sum_{|j| \leq 2} \int_t^T \E{\E{|D^j_x \psi_s(x', z')|^p}_{(x', z') = ( \Xetx, \Zetz )}} ds.
\end{align*}

\begin{proof}

$R_t(x, z)$ solves the following BSPDE
\begin{align}
\label{equation: BSPDE, remainder dual process}
-dR_t 
&= (\mG^\epsilon R_t + \mG_S \psi_t) dt + \langle (\psi_t +  R_t) h , d \B \rangle + \langle \alpha \sigma^* \nabla_x (\psi_t + R_t) , d \B \rangle,
&\qquad R_T = 0.
\end{align}
Existence of the solution $R$ and its derivatives as well as the polynomial growth all follow from Lemma \ref{lemma: Rozovskii SPDE result}. 
The parabolic condition of Lemma \ref{lemma: Rozovskii SPDE result} holds because $I - \alpha^* \alpha \succeq 0$, where $I$ is the identity matrix. 
From Lemma \ref{lemma: BDSDE implies BSPDE solution}, the solution, $R_t(x, z)$, has a representation in terms of a FBDSDE, $R_t(x, z) = \theta_t^{t, x, z}$. 
Where $\theta_t$ is the first component of the tuple of processes $(\theta_s, \gamma_s^\eps, \eta_s)$ satisfying the BDSDE
\EquationAligned{
\label{equation: BDSDE, remainder process}
-d \theta_s^{t, x, z} 
&= \mG_S \psi_s (\Xetx, \Zetz) ds \\
&+ \langle \psi_s h (\Xetx , \Zetz) , d \B[s] \rangle 
+ \langle \theta^{t, x, z}_s h (\Xetx , \Zetz) , d \B[s] \rangle \\
&+ \langle \alpha \sigma^* \nabla_x \psi_s (\Xetx, \Zetz) , d \B[s] \rangle 
+ \langle \alpha \eta_s^{t, x, z} , d \B[s] \rangle \\
&-  \langle \eta_s^{t, x, z} , dW_s \rangle
- \langle \gamma_s^{\epsilon; t, x, z} , dV_s \rangle , \\
\theta_T^{t, x, z} &= 0 ,
}
and $(\Xetx, \Zetz)$ is a joint diffusion process satisfying the SDEs
\EquationAligned*{
d\Xetx 
&= b(\Xetx, \Zetz) ds + \sigma(\Xetx, \Zetz) dW_s, \qquad s \geq t , \\
\Xetx &= x, \qquad s \leq t , \\[2ex]
d\Zetz 
&= \frac{1}{\epsilon^2} f(\Xetx, \Zetz) ds + \frac{1}{\epsilon} g(\Xetx, \Zetz) dV_s, \qquad s \geq t , \\
\Zetz &= z, \qquad s \leq t,
}
where we choose $(W, V)$ and $B$ to be independent standard Brownian motions. 
This is necessary when working with a stochastic representation of Eq. \ref{equation: BSPDE, remainder dual process}. 
The second and third components of the tuple $(\theta^{t, x, z}_t, \gamma^{\epsilon; (t, x, z)}_t, \eta^{t, x, z}_t)$, have representations as 
\begin{align*}
\gamma^{\epsilon; (t, x, z)}_t = \frac{1}{ \epsilon } g^* \nabla_z R_t (x, z) 
\qquad 
\text{and}
\qquad 
\eta^{t, x, z}_t = \sigma^* \nabla_x R_t (x, z).  
\end{align*}
Let $\mathcal{A}$ be the integrand for the backward stochastic integral in Eq. \ref{equation: BDSDE, remainder process}; it takes the following definition 
\begin{align*}
\mathcal{A} 
= \psi_s h (\Xetx , \Zetz)
+ \theta^{t, x, z}_s h (\Xetx , \Zetz) 
+ \alpha \sigma^* \nabla_x \psi_s (\Xetx, \Zetz)
+ \alpha \eta_s^{t, x, z} ,
\end{align*}
or stripping function arguments and superscripts, 
\begin{align*}
\mathcal{A} 
= \psi_s h
+ \theta_s h
+ \alpha \sigma^* \nabla_x \psi_s
+ \alpha \eta_s .
\end{align*}
We now consider the $p$-th moment of $\theta_t$, 
\EquationAligned{
\label{equation: p-th moment of error - Ito approach}
\E{| \theta_t |^p  }
&= \int_t^T \E{ p |\theta_s|^{p-2} \theta_s \mG_S \psi_s } ds
+ \frac{p (p - 1)}{2} \int_t^T \E{ | \theta_s |^{p - 2} | \mathcal{A} |^2 } ds \\
&- \frac{p (p - 1)}{2} \int_t^T \E{ | \theta_s |^{p - 2} | \eta_s |^2 } ds 
- \frac{p (p - 1)}{2} \int_t^T \E{ | \theta_s |^{p - 2} | \gamma_s |^2 } ds.
}

Using the fact that $\theta, \psi$ are real-valued functions, and $b, \sigma \in C^{k, l}_b$, applying Young's inequality to the first term on the right side of Eq. \ref{equation: p-th moment of error - Ito approach} yields, 
\begin{align*}
\int_t^T \E{ p |\theta_s|^{p-2} \theta_s \mG_S \psi_s } ds 
\leq \frac{p}{2} \int_t^T \E{ |\theta_s|^p } ds  + \frac{p}{2} \int_t^T \E{ |\theta_s|^{p-2} | \mG_S \psi_s |^2 } ds.
\end{align*}
Application of H\"older's inequality and Young's inequality to the last term gives,
\EquationAligned{
\label{equation: remainder term, time integral, holder application}
\frac{p}{2} \int_t^T \E{ |\theta_s|^{p-2} | \mG_S \psi_s |^2 } ds
&\leq \frac{p - 2}{2} \int_t^T \E{ |\theta_s|^p } ds + \int_t^T \E{ |\mG_S \psi_s|^p } ds. 
}
Application of H\"older's inequality and the use of the boundedness of $b$ and $a$, then the tower property of conditional expectation, and the Markov property of $(\Xe, \Ze)$ gives the following bound for the last term in Eq. \ref{equation: remainder term, time integral, holder application}, 
\begin{align*}
\int_t^T \E{ |\mG_S \psi_s|^p } ds 
&\lesssim \int_t^T \sum_{|j| \leq 2} \E{ \left| D^j_x \psi(\Xetx, \Zetz) \right|^p } ds \\ 
&= \int_t^T \sum_{|j| \leq 2} \E{ \E{ \left| D^j_x \psi(\Xetx, \Zetz) \right|^p }[\mF^W_s \vee \mF^V_s] } ds \\
&= \int_t^T \sum_{|j| \leq 2} \E{ \E{ \left| D^j_x \psi(x', z') \right|^p }_{(x', z') = (\Xetx, \Zetz)} } ds .
\end{align*}
Therefore the first term on the right side of Eq. \ref{equation: p-th moment of error - Ito approach} is bounded by, 
\EquationAligned{
\label{equation: remainder term, time integral, final estimate}
\int_t^T \E{ p |\theta_s|^{p-2} \theta_s \mG_S \psi_s } ds 
&\lesssim (p - 1) \int_t^T \E{ |\theta_s|^p } ds  \\
&+ \int_t^T \sum_{|j| \leq 2} \E{ \E{ \left| D^j_x \psi(x', z') \right|^p }_{(x', z') = (\Xetx, \Zetz)} } ds. 
}

Now addressing the second term on the right side of Eq. \ref{equation: p-th moment of error - Ito approach}, expanding the inner product $| \mathcal{A} |^2 = \langle \mathcal{A} , \mathcal{A} \rangle$ and separating terms using Young's inequality  with values $\lambda_1, \hdots, \lambda_6 > 0$ to be chosen later, we get 
\EquationAligned*{
\langle \mathcal{A} , \mathcal{A} \rangle
&\leq \left(1 + \frac{1}{\lambda_1} + \frac{1}{\lambda_2} + \frac{1}{\lambda_3} \right) | \psi_s h |^2
+ \left(1 + \lambda_1 + \frac{1}{\lambda_4} + \frac{1}{\lambda_5} \right) | \theta_s h |^2 \\
&+ \left(1 + \lambda_2 + \lambda_4 + \frac{1}{\lambda_6} \right) | \alpha \sigma^* \nabla_x \psi_s |^2
+ \left(1 + \lambda_3 + \lambda_5 + \lambda_6 \right) | \alpha \eta_s |^2,
}
Therefore the second term on the right side of Eq. \ref{equation: p-th moment of error - Ito approach} is bounded by 
\begin{align}
\frac{p (p - 1)}{2} \int_t^T \E{ | \theta_s |^{p - 2} | \mathcal{A} |^2 } ds
&\leq \left(1 + \frac{1}{\lambda_1} + \frac{1}{\lambda_2} + \frac{1}{\lambda_3} \right) \frac{p (p - 1)}{2} \int_t^T \E{ | \theta_s |^{p - 2} | \psi_s h |^2 } ds \label{equation: remainder term, Young's expansion, 1st term} \\
&+ \left(1 + \lambda_1 + \frac{1}{\lambda_4} + \frac{1}{\lambda_5} \right) \frac{p (p - 1)}{2} \int_t^T \E{ | \theta_s |^{p - 2} | \theta_s h |^2 } ds \label{equation: remainder term, Young's expansion, 2nd term} \\
&+ \left(1 + \lambda_2 + \lambda_4 + \frac{1}{\lambda_6} \right) \frac{p (p - 1)}{2} \int_t^T \E{ | \theta_s |^{p - 2} | \alpha \sigma^* \nabla_x \psi_s |^2 } ds \label{equation: remainder term, Young's expansion, 3rd term} \\
&+ \left(1 + \lambda_3 + \lambda_5 + \lambda_6 \right) \frac{p (p - 1)}{2} \int_t^T \E{ | \theta_s |^{p - 2} | \alpha \eta_s |^2 } ds . \label{equation: remainder term, Young's expansion, 4th term}
\end{align}

We now consider pairing the term given by Eq. \ref{equation: remainder term, Young's expansion, 4th term} and the third term on the right side of Eq. \ref{equation: p-th moment of error - Ito approach},
\EquationAligned{
\label{equation: error term, dual identity in Ito expansion}
\underbrace{\left(1 + \lambda_3 + \lambda_5 + \lambda_6 \right)}_{\equiv \Lambda} \frac{p (p - 1)}{2} \int_t^T \E{ |\theta_s|^{p - 2} | \alpha \eta_s |^2 } ds  
- \frac{p (p - 1)}{2} \int_t^T \E{ |\theta_s|^{p - 2} | \eta_s |^2 } ds \\
= \frac{p (p - 1)}{2} \int_t^T \E{ |\theta_s|^{p - 2} \left( \Lambda \eta_s^* \alpha^* \alpha \eta_s - \eta_s^* \op{Id} \eta_s \right) } ds \\ 
= \frac{p (p - 1)}{2} \int_t^T \E{ |\theta_s|^{p - 2} \left( \eta_s^* \left( \Lambda \alpha^* \alpha - \op{Id} \right) \eta_s \right) } ds.
}
The constant matrix $\alpha^* \alpha - \op{Id}$ is negative definite and we can choose $\lambda_3, \lambda_5, \lambda_6 > 0$, small enough such that $\Lambda \alpha^* \alpha - \op{Id} \prec 0$.

Turning our attention to the three terms of Eqs. \ref{equation: remainder term, Young's expansion, 1st term}, \ref{equation: remainder term, Young's expansion, 2nd term}, and \ref{equation: remainder term, Young's expansion, 3rd term}, we use the same technique as in Eq. \ref{equation: remainder term, time integral, holder application} with $|h|_\infty < \infty$ on the first term (Eq. \ref{equation: remainder term, Young's expansion, 1st term}) to get
\begin{multline}
\label{equation: remainder term, Young's expansion, 1st term, final estimate}
\left(1 + \frac{1}{\lambda_1} + \frac{1}{\lambda_2} + \frac{1}{\lambda_3} \right) \frac{p (p - 1)}{2} \int_t^T \E{ |\theta_s|^{p - 2} | \psi_s h |^2  } ds  
\lesssim \int_t^T \E{ | \theta_s |^{p - 2} | \psi_s h |^2 } ds \\
\lesssim \int_t^T \E{ | \theta_s |^p } ds + \int_t^T \sum_{|j| \leq 0} \E{ \E{ \left| D^j_x \psi(x', z') \right|^p }_{(x', z') = (\Xetx, \Zetz)} } ds.
\end{multline}
Similarly, again using $| h |_\infty < \infty$, Eq. \ref{equation: remainder term, Young's expansion, 2nd term} is bounded by
\begin{align}
\label{equation: remainder term, Young's expansion, 2nd term, final estimate}
\left(1 + \lambda_1 + \frac{1}{\lambda_4} + \frac{1}{\lambda_5} \right) \frac{p (p - 1)}{2} \int_t^T \E{ |\theta_s|^{p - 2} | \theta_s h |^2  } ds  
\lesssim \int_t^T \E{ | \theta_s |^p } ds.
\end{align}
And now Eq. \ref{equation: remainder term, Young's expansion, 3rd term} using $| \sigma |_\infty < \infty$, 
\begin{multline}
\label{equation: remainder term, Young's expansion, 3rd term, final estimate}
\left(1 + \lambda_2 + \lambda_4 + \frac{1}{\lambda_6} \right) \frac{p (p - 1)}{2} \int_t^T \E{ |\theta_s|^{p - 2} | \alpha \sigma^* \nabla_x \psi_s |^2 } ds   \\
\lesssim \int_t^T \E{ | \theta_s |^p } ds + \int_t^T \sum_{|j| \leq 1} \E{ \E{ \left| D^j_x \psi(x', z') \right|^p }_{(x', z') = (\Xetx, \Zetz)} } ds.
\end{multline}

Collecting the bounds of Eqs. \ref{equation: remainder term, time integral, final estimate}, \ref{equation: error term, dual identity in Ito expansion}, \ref{equation: remainder term, Young's expansion, 1st term, final estimate}, \ref{equation: remainder term, Young's expansion, 2nd term, final estimate}, and \ref{equation: remainder term, Young's expansion, 3rd term, final estimate} we get for Eq. \ref{equation: p-th moment of error - Ito approach}, 
\EquationAligned*{
\E{ | \theta_t |^p } 
\lesssim \int_t^T \E{ | \theta_s |^p } ds 
&+ \int_t^T \sum_{|j| \leq 2} \E{ \E{ \left| D^j_x \psi(x', z') \right|^p }_{(x', z') = (\Xetx, \Zetz)} } ds \\
&+ \frac{p (p - 1)}{2} \int_t^T \E{ |\theta_s|^{p - 2} \left( \eta_s^* \left( \Lambda \alpha^* \alpha - \op{Id} \right) \eta_s \right) } ds \\
&- \frac{p (p - 1)}{2} \int_t^T \E{ | \theta_s |^{p - 2} | \gamma_s |^2 } ds.
}
Rearranging this equation gives
\EquationAligned*{
\E{ | \theta_t |^p } 
- \frac{p (p - 1)}{2} \int_t^T \E{ |\theta_s|^{p - 2} \left( \eta_s^* \left( \Lambda \alpha^* \alpha - \op{Id} \right) \eta_s \right) } ds 
+ \frac{p (p - 1)}{2} \int_t^T \E{ | \theta_s |^{p - 2} | \gamma_s |^2 } ds \\
\lesssim \int_t^T \E{ | \theta_s |^p } ds 
+ \int_t^T \sum_{|j| \leq 2} \E{ \E{ \left| D^j_x \psi(x', z') \right|^p }_{(x', z') = (\Xetx, \Zetz)} } ds.
}
From the fact that $\Lambda \alpha^* \alpha - \op{Id} \prec 0$, the subtraction of the second term on the left side of the equation is a non-negative value. 
The third term on the left side of the equation is also non-negative, and therefore we can drop them from the inequality to get
\EquationAligned*{
\E{ | \theta_t |^p } 
\lesssim \int_t^T \E{ | \theta_s |^p } ds 
+ \int_t^T \sum_{|j| \leq 2} \E{ \E{ \left| D^j_x \psi(x', z') \right|^p }_{(x', z') = (\Xetx, \Zetz)} } ds,
}
Now applying Gr\"onwall's lemma yields, 
\EquationAligned*{
\E{ | \theta_t |^p } 
&\lesssim \int_t^T \sum_{|j| \leq 2} \E{ \E{ \left| D^j_x \psi(x', z') \right|^p }_{(x', z') = (\Xetx, \Zetz)} } ds. 
}
Using the fact that the solution to the BDSDE provides the classical solution to the BSPDE, $R_t(x, z) = \theta_t^{t, x, z}$, we get the desired result. 
\end{proof}
}
}

\subsection{Estimates of Dual and Filter Error}
\label{subsection: estimates of dual and filter error}

We now complete the final estimates that lead to the proof of Theorem \ref{theorem: main result of quantitative convergence}.
The remaining lemmas parallel those in \cite{Imkeller:2013}; we provide the statements, which have slightly different assumptions, but otherwise refer to \cite{Imkeller:2013} for the proofs.

{
\lem{
Assume \ref{assumption: positive recurrence}, \ref{assumption: uniform ellipticity}, \hyperlink{HF}{$HF^{8, 4}$}, $b \in C^{7, 4}_b$, $\sigma \in C^{8, 4}_b$, \hyperlink{HO}{$HO^{8, 4}$}, and $\varphi \in C^7_b(\R{m}; \R{})$.
Then for every $p \geq 1$ there exists $q > 0$, such that 
\begin{align*}
\sup_{0 \leq t \leq T} \E[\Pe]{ \left| v^{\epsilon, T, \varphi}_t(x, z) - v^{0, T, \varphi}_t(x) \right|^p }
\lesssim \epsilon^p ( 1 + |x|^q + |z|^q ) \norm{\varphi}{4, \infty}^p .
\end{align*}

\begin{proof}
First we collect the conditions in reverse order of our main estimates. 
\begin{enumerate}[(i)]
\item For the solution of $R$ in Lemma \ref{lemma: quantitative remainder term}, we require \hyperlink{HF}{$HF^{3, 3}$}, \hyperlink{HS}{$HS^{3, 3}$}, $\sigma \in C^{4, 4}_b$, \hyperlink{HO}{$HO^{4, 4}$}, and $\psi \in C^{0, 5, 3}$. 
The polynomial growth will be satisfied. 
\item For the solution of $\psi \in C^{0, 5, 3}$ in Lemma \ref{lemma: quantitative corrector term}, we require \hyperlink{HF}{$HF^{6, 4}$}, \hyperlink{HS}{$HS^{6, 4}$}, \hyperlink{HO}{$HO^{6, 4}$}, and $v^0 \in C^{0, 6}$. 
The conditions on $\overbar{\sigma}, \overbar{a}, \overbar{b}$, and $\overbar{h}$ will already be covered by the stronger conditions just stated. 
The polynomial growth will also be satisfied. 
And we require \ref{assumption: positive recurrence} and \ref{assumption: uniform ellipticity}.
\item For the solution of $v^0 \in C^{0, 6}$ in Lemma \ref{lemma: quantitative averaged term}, we require $\overbar{b}, \overbar{a} \in C^7_b$, that $\overbar{h}, \overbar{\sigma} \in C^8_b$, and $\varphi \in C^7_b$.
\item Using Lemma \ref{lemma: regularity and integrability of simple averaged coefficients (chapter: quantitative)}, for $\overbar{h}, \overbar{\sigma} \in C^8_b$ requires \hyperlink{HF}{$HF^{8, 3}$}, \hyperlink{HO}{$HO^{8, 1}$}, and $\sigma \in C^{8, 1}_b$.
And this also implies $\overbar{a} \in C^7_b$. 
For $\overbar{b} \in C^7_b$, we need $b \in C^{7, 1}_b$. 
We also require \ref{assumption: positive recurrence} and \ref{assumption: uniform ellipticity}.
\item Therefore the sufficient conditions are \ref{assumption: positive recurrence}, \ref{assumption: uniform ellipticity}, \hyperlink{HF}{$HF^{8, 4}$}, $b \in C^{7, 4}_b$, $\sigma \in C^{8, 4}_b$, \hyperlink{HO}{$HO^{8, 4}$}, and $\varphi \in C^7_b$. 
\end{enumerate}
For the remainder of the proof, see \cite[Lemma 6.4, p.2318]{Imkeller:2013}.
\end{proof}
}
}

We now show that the moment estimate of the difference of $v^\epsilon$ and $v^0$ continues to hold under the original measure $\Q{}$. 
{
\lem{
\label{lemma: dual error moment bound, under original measure}
Assume \ref{assumption: positive recurrence}, \ref{assumption: uniform ellipticity}, \hyperlink{HF}{$HF^{8, 4}$}, $b \in C^{7, 4}_b$, $\sigma \in C^{8, 4}_b$, \hyperlink{HO}{$HO^{8, 4}$}, and $\varphi \in C^7_b(\R{m}; \R{})$.
Then for every $p \geq 1$ there exists $q > 0$, such that 
\begin{align*}
\sup_{0 \leq t \leq T} \E[\Q{}]{ \left| v^{\epsilon, T, \varphi}_t(x, z) - v^{0, T, \varphi}_t(x) \right|^p }
\lesssim \epsilon^p ( 1 + |x|^q + |z|^q ) \norm{\varphi}{4, \infty}^p .
\end{align*}

\begin{proof}
Using Lemma \ref{lemma: moments of Girsanov (chapter: quantitative)}, the proof follows \cite[Lemma 6.5, p.2319]{Imkeller:2013}.
\end{proof}
}
}

{
\lem{
\label{lemma: moment estimate of unnormalized measure error}
Assume \ref{assumption: positive recurrence}, \ref{assumption: uniform ellipticity}, \hyperlink{HF}{$HF^{8, 4}$}, $b \in C^{7, 4}_b$, $\sigma \in C^{8, 4}_b$, \hyperlink{HO}{$HO^{8, 4}$}, and $\varphi \in C^7_b(\R{m}; \R{})$.
Additionally, assume that the initial distribution $\Q{}_{(\Xe_0, \Ze_0)}$ has finite moments of every order. 
Then for every $p \geq 1$ there exists $q > 0$, such that 
\begin{align*}
\E[\Q{}]{ \left| \rho^{\epsilon, x}_T (\varphi) - \rho^0_T (\varphi) \right|^p }
\lesssim \epsilon^p \norm{\varphi}{4, \infty}^p .
\end{align*}

\begin{proof}
The proof is the same as \cite[Lemma 6.6, p.2320]{Imkeller:2013}
\end{proof}
}
}

{
\lem{
\label{lemma: inverse moment of unnormalized measure on unit function} 
Let $p \geq 1$ and assume $h$ is bounded. 
Then 
\begin{align*}
\sup_{\epsilon \in (0, 1]} \sup_{0 \leq t \leq T} \left( \E[\Q{}]{ | \rho^{\epsilon, x}_t (1) |^{-p} } + \E[\Q{}]{ | \rho^0_t(1) |^{-p} } \right) 
< \infty.
\end{align*}
\begin{proof}
For the first term, the proof is the same as \cite[Lemma 6.7, p.2321]{Imkeller:2013}.  
The same is true for the second term using the definitions given in Lemma \ref{lemma: existence of averaged unnormalized filter}. 
\end{proof}
}
}

{
\lem{
\label{lemma: moment estimate of normalized measure error}
Assume \ref{assumption: positive recurrence}, \ref{assumption: uniform ellipticity}, \hyperlink{HF}{$HF^{8, 4}$}, $b \in C^{7, 4}_b$, $\sigma \in C^{8, 4}_b$, \hyperlink{HO}{$HO^{8, 4}$}, and $\varphi \in C^7_b(\R{m}; \R{})$.
Additionally, assume that the initial distribution $\Q{}_{(\Xe_0, \Ze_0)}$ has finite moments of every order. 
Then for every $p \geq 1$ there exists $q > 0$, such that 
\begin{align*}
\E[\Q{}]{ \left| \pi^{\epsilon, x}_T (\varphi) - \pi^0_T (\varphi) \right|^p }
\lesssim \epsilon^p \norm{\varphi}{4, \infty}^p .
\end{align*}

\begin{proof}
The proof is the same as \cite[Lemma 6.8, p.2321]{Imkeller:2013}
\end{proof}
}
}

Observing that the bound in the result of Lemma \ref{lemma: moment estimate of normalized measure error} only depends on $\norm{\varphi}{4, \infty}^p$, even though the assumption requires $\varphi \in C^7_b$, encourages us to instead approximate a fixed test function $\varphi \in C^4_b$ by a sequence $(\varphi^n \in C^7_b)$ in the $\norm{\cdot}{4, \infty}$-norm, and take advantage of the fact that $\pi^{\epsilon, x}_T$ and $\pi^0_T$ are $\Q{}$-a.s. equal to probability measures. 
Therefore we can relax this condition in Lemma \ref{lemma: moment estimate of normalized measure error} slightly with the following corollary. 

{
\cor{
\label{corollary: moment estimate of normalized measure error}
Assume \ref{assumption: positive recurrence}, \ref{assumption: uniform ellipticity}, \hyperlink{HF}{$HF^{8, 4}$}, $b \in C^{7, 4}_b$, $\sigma \in C^{8, 4}_b$, and \hyperlink{HO}{$HO^{8, 4}$}.
Additionally, assume that the initial distribution $\Q{}_{(\Xe_0, \Ze_0)}$ has finite moments of every order. 
Then for any $p \geq 1$ we have that for every $\varphi \in C^4_b(\R{m}; \R{})$,
\begin{align*}
\E[\Q{}]{ \left| \pi^{\epsilon, x}_T (\varphi) - \pi^0_T (\varphi) \right|^p }
\lesssim \epsilon^p \norm{\varphi}{4, \infty}^p .
\end{align*}
}
}

The next lemma shows that indeed we have weak convergence of $\pi^{\epsilon, x}$ to $\pi^0$.

{
\lem{
\label{lemma: convergence of normalized measure error under metric generating weak topology}
Assume \ref{assumption: positive recurrence}, \ref{assumption: uniform ellipticity}, \hyperlink{HF}{$HF^{8, 4}$}, $b \in C^{7, 4}_b$, $\sigma \in C^{8, 4}_b$, and \hyperlink{HO}{$HO^{8, 4}$}.
Additionally, assume that the initial distribution $\Q{}_{(\Xe_0, \Ze_0)}$ has finite moments of every order. 
Then there exists a metric $d$ on the space of probability measures on $\R{m}$ that generates the topology of weak convergence, such that 
\begin{align*}
\E[\Q{}]{ d( \pi^{\epsilon, x}_T , \pi^0_T ) }
\lesssim \epsilon .
\end{align*}

\begin{proof}
To achieve this result, we borrow the argument following \cite[Corollary 6.9, p.2322]{Imkeller:2013}. 
\end{proof}
}
}

\section*{Acknowledgement}
R.B. and N.S.N. acknowledge partial support for this work from the Air Force Office of Scientific Research under grant number FA9550-17-1-0001, and N.S.N. acknowledges partial support from the National Sciences and Engineering Research Council Discovery grant 50503-10802.

\sloppy{
\printbibliography[heading=bibintoc]
}
\end{document}